\documentclass[11pt]{amsart} 
\usepackage{curves}
\usepackage{euler}
\usepackage{euscript}
\usepackage{multicol} 
\usepackage{amssymb}
\usepackage{amsmath} 
\usepackage{times} 

\newtheorem{theorem}{Theorem}[section] 
\newtheorem{lemma}[theorem]{Lemma}
\newtheorem{proposition}[theorem]{Proposition}

\newtheorem{corollary}[theorem]{Corollary}
\numberwithin{equation}{section}

\begin{document}
\title[RH approach to Toeplitz operators and orthogonal polynomials]{A Riemann-Hilbert 
approach to some theorems on Toeplitz operators and orthogonal polynomials}
\author[P. Deift, J. \"{O}stensson]{Percy Deift and J\"{o}rgen \"{O}stensson\\
\\
\textit{This paper is dedicated to Barry Simon on the occasion of his 60th birthday in appreciation for all 
that he has taught us.}}
\begin{abstract}
In this paper the authors show how to use Riemann-Hilbert techniques to prove various results, some 
old, some new, in the theory of Toeplitz operators and orthogonal polynomials on the unit circle (OPUC's).
There are four main results: the first concerns the approximation of the inverse of a Toeplitz operator 
by the inverses of its finite truncations. The second concerns a new proof of the `hard' part of Baxter's 
theorem, and the third concerns the Born approximation for a scattering problem on the lattice $\mathbb{Z}_+$.
The fourth and final result concerns a basic proposition of Golinskii-Ibragimov arising in their analysis of the
Strong Szeg\"{o} Limit Theorem.
\end{abstract} 
\maketitle
\setcounter{section}{-1}
\section{Introduction.} 
Let $d\mu$ be a probability measure on the unit circle $\Gamma = \{z \in \mathbb{C} : |z| = 1\}$ and let 
$\Phi_n = z^n + ...\,$, $n \geq 0$, be the (monic) orthogonal polynomials (OPUC's) associated with $d\mu$, 
$\int_{\Gamma} \Phi_m(z)\,\overline{\Phi_n(z)}\,d\mu = 0$, $m \not= n$, $m, n \geq 0$ (see \cite{Sz}). 
Let $\alpha = (\alpha_n)_{n \in \mathbb{Z}_+}$ denote the vector
of \textit{Verblunsky coefficients} $\alpha_n = - \overline{\Phi_{n+1}(0)}$, $n \geq 0$. By Verblunsky's theorem (see \cite{Sim1}), the map 
$V : d\mu \mapsto \alpha$ is a bijection from the probability measures on $\Gamma$ onto $\times_{j =0}^{\infty}\,\mathbb{D}$, 
where $\mathbb{D} = \{z \in \mathbb{C} : |z| < 1\}$ is the (open) unit disc in $\mathbb{C}$. Following Cantero, Moral and Vel\'{a}zquez \cite{CMV}, 
we may, given $\alpha$, construct a (pentadiagonal) unitary matrix operator $U = U(\alpha)$ in $l^2_+ = l^2(\mathbb{Z}_+)$ (the so-called CMV matrix) 
with the following property: $e_0 = (1, 0, ...)^{T}$ is a cyclic vector for $U$, i.e. $\overline{<U^k\,e_0>}_{-\infty < k < \infty} = l^2_+$, and 
the associated spectral measure for $U$ is precisely $d\mu = V^{-1}(\alpha)$. With this construction, Verblunsky's theorem becomes a 
result in spectral/inverse spectral theory: Indeed, let $\mathcal{S}$ denote the map from CMV matrices $U$ to their spectral 
measures $d\mu$ on $\Gamma$,
\begin{equation}
U \longmapsto d\mu 
\end{equation}
and let $\mathcal{I}$ denote the map from measures $d\mu$ on $\Gamma$ to their associated CMV matrices $U = U(V(d\mu))$,
\begin{equation}
d\mu \longmapsto U(V(d\mu)).
\end{equation}
Then $\mathcal{S}$ and $\mathcal{I}$ are inverse to each other. The above correspondence, which is the analog 
for the unit circle of the well-known correspondence between measures on the line and Jacobi operators (see e.g. \cite{D2}), divides the study
of OPUC's naturally into two parts: the \textit{direct problem} (equivalently, the study of the properties of $\mathcal{S}$) and the 
\textit{inverse problem} (equivalently, the study of the properties of $\mathcal{I}$). This is the approach taken in Simon's new book \cite{Sim1,Sim2}:
Part 1 focuses on $\mathcal{I}$ and Part 2 focuses on $\mathcal{S}$. The goal of the present paper is to show that the study of the map 
$\mathcal{I}$ is greatly facilitated by using Riemann-Hilbert (RH) techniques. We will do this by producing new and transparent RH proofs of some
classical and central theorems in the subject: En route, we will also derive some new results.  

Denote by $H_{\pm}$ the closed subspaces of $L^2\left(\Gamma\right)$ 
consisting of functions $u$ whose negative/non-negative Fourier coeffients are zero, and let $P_{\pm}\,:\, L^2\left(\Gamma\right) \rightarrow H_{\pm}$ 
be the associated orthogonal projections.
Given a function $\varphi \in L^{\infty}\left(\Gamma\right)$ we define the associated \textit{Toeplitz operator with symbol $\varphi$},
$T\left(\varphi\right)\,:\, H_+ \rightarrow H_+$, by the formula 
\begin{equation}
T(\varphi)\,u = P_+(\varphi\,u), \quad u \in H_+.
\end{equation}
In terms of the Fourier coefficients $\varphi_k =  \widehat{\varphi}\,(k) = \int_{-\pi}^{\pi} e^{-ik\theta}\,\varphi(e^{i\theta})\frac{d\theta}{2\pi}$ 
the Toeplitz operator becomes a truncated discrete convolution:
\begin{equation}
T(\varphi)\,z^k = \sum_{j=0}^{\infty} \varphi_{j-k}\,z^j, \quad z \in \Gamma, \quad k \in \mathbb{Z}_+.
\end{equation}
Let $T(\varphi)_{jk} = \varphi_{j-k}$. Then the \textit{Toeplitz matrix} $\left(T(\varphi)_{jk}\right)_{j,k = 0}^{\infty} = 
\left(\varphi_{j-k}\right)_{j,k = 0}^{\infty}$ is the matrix representation of $T(\varphi)$ in the standard basis 
$\left(z^k\right)_{k=0}^{\infty}$ for $H_+$.
For $n \geq 0$, let $\mathcal{P}_n = \left\{\sum_{j=0}^n a_j\,z^j \right\}$ denote the subspace of $L^2\left(\Gamma\right)$ consisting 
of polynomials of degree less than or equal to $n$, and $P_n\,:\,L^2\left(\Gamma\right) \rightarrow \mathcal{P}_n$ the corresponding orthogonal 
projection. Define the \textit{n'th truncation} of the Toeplitz operator $T(\varphi)$ to be the map 
$T_n = T_n(\varphi) = P_n\,T(\varphi)|_{\mathcal{P}_n}$.

In the following we will be interested only in symbols $\varphi$ belonging to the so-called Beurling class $W_{\nu}$ (compare \cite{Sim1}). 
The basic definitions are as follows. We call a sequence $\nu = \left(\nu_k\right)_{k \in \mathbb{Z}}$ a \textit{Beurling weight} if 
it has the properties:\\
$$
\begin{array}{lll}
(i) \quad &\nu_{j} \geq 1, \quad &j \in \mathbb{Z}\\ 
(ii) \quad &\nu_{j} = \nu_{-j}, \quad &j \in \mathbb{Z}\\
(iii) \quad &\nu_{j+k} \leq \nu_{j}\,\nu_{k}, \quad &j,k \in \mathbb{Z}\\
\end{array}
$$
The Beurling class is defined as 
\begin{equation*}
W_{\nu} = \bigg\{\varphi \in L^1\left(\Gamma\right) \,:\, \sum_{j \in \mathbb{Z}} \nu_j \,|\varphi_j| < \infty \bigg\}.
\end{equation*}
By standard subadditivity arguments it follows that 
\begin{equation}
A(\nu) = \lim_{k \rightarrow \infty} \frac{\log \nu_k}{k} = \inf_{k \in \mathbb{N}} \frac{\log \nu_{k}}{k} 
\end{equation} 
exists. Note, in particular, that $A(\nu) \geq 0$ and also that $\nu_k \geq e^{|k|\,A(\nu)}, k \in \mathbb{Z}$. In case $A(\nu) = 0$, we say 
that $\nu$ is a \textit{strong} Beurling weight. It is easy to see that $W_{\nu}$ becomes a Banach algebra if equipped with the norm
\begin{equation}
  \label{Banachnorm}
||\varphi||_{\nu} = \sum_{j \in \mathbb{Z}} \nu_j \,|\varphi_j|. 
\end{equation}
Canonical examples are given by the exponential weights $\nu_j = \gamma^{|j|}, \gamma \geq 1$, and the algebra $W^{\alpha}$ associated 
with (strong) Beurling weight $\nu_j = (1 + |j|)^{\alpha}$, $\alpha \geq 0$.  The space $W^{0}$ is the \textit{standard Wiener algebra}. 
Note that $W_{\nu} \subset W^{0}$ for any Beurling weight $\nu$.

It is a well-known theorem, due to Krein, that if $\varphi \in W^{0}$, then $T(\varphi)$ is invertible if and only if
$\varphi(z) \ne 0$ for all $z \in \Gamma$ and $wind(\varphi, 0) = 0$. 
In this case, the inverse is given by
\begin{equation}
T(\varphi)^{-1} = T\left(\frac{1}{\varphi_+}\right)\,T\left(\frac{1}{\varphi_-}\right),
\end{equation}    
where $\varphi = \varphi_+\,\varphi_-$ is the Wiener-Hopf factorization of $\varphi$, i.e.
$\varphi_+$ extends to a non-vanishing function analytic in the interior of the unit circle and $\varphi_-$ to a non-vanishing 
function, with $\varphi_-(\infty) = 1$, analytic in the exterior of the unit circle. Said differently,
\begin{equation*}
m(z) = 
\left\{
\begin{array}{ll}
\varphi_+(z)&, |z| < 1,\\
\varphi_-^{-1}(z)&, |z| > 1,   
\end{array}
\right.
\end{equation*}
is the solution of the (scalar) Riemann-Hilbert Problem (RHP) $(\Gamma, v = \varphi)$ (see below). It is not difficult to see that, under the above 
conditions on $\varphi$, such a factorization exists and that the extensions are uniquely given by $\varphi_{\pm} = \exp\,\{\pm C(\log \varphi)\}$.

Suppose that $\varphi \in W_{\nu}$. Let us denote by $\mathcal{R}_{\nu}$ the annulus
\begin{equation*}
\mathcal{R}_{\nu} = \left\{z \in \mathbb{C} : e^{-A(\nu)} \leq |z| \leq e^{A(\nu)}\right\}.
\end{equation*}
It is then easy to see that $\varphi$ extends to a function analytic in the interior of $\mathcal{R}_{\nu}$ and continuous up to the boundary.  
Using basic facts from the Gelfand theory of commutative Banach algebras one can prove that the spectrum $\sigma(\varphi)$ of $\varphi$ 
equals $\varphi(\mathcal{R}_{\nu})$, i.e. if $\varphi(z) \not= 0$ for $z \in \mathcal{R}_{\nu}$, then $\varphi^{-1} \in W_{\nu}$. 
Furthermore, if in addition to the assumption that $\varphi \in W_{\nu}$ is non-vanishing on $\mathcal{R}_{\nu}$ we impose the condition 
that $wind(\varphi,0) = 0$, then $\log \varphi \in W_{\nu}$. This follows from the following basic fact, see \cite{Doug}: 
Let us denote by $G\mathcal{B}$ the group of invertible elements of a commutative Banach algebra $\mathcal{B}$ and 
by $G_{0}\mathcal{B}$ the (connected) component in $G\mathcal{B}$ containing the identity. Then, $G_{0}\mathcal{B}$ coincides 
with $\exp \mathcal{B}$. Indeed, write $\varphi(z) = \sum_{j \in \mathbb{Z}} a_j\,z^j$ and introduce the sequence of rational 
approximations $\varphi^{(N)}(z) = \sum_{j = -N}^{N} a_j \,z^j$. Clearly then $\varphi^{(N)} \in W_{\nu}$, and 
$\varphi^{(N)} \rightarrow \varphi$ in $W_{\nu}$. It follows that, for $N$ sufficiently large, $\varphi^{(N)}$ is non-vanishing on 
$\mathcal{R}_{\nu}$ with $wind(\varphi^{(N)},0) = 0$. Clearly then, for such $N$,  
\begin{equation*}
\varphi^{(N)}(z) = c\,\frac{\Pi_{j = 1}^{N}(z-\alpha_j)\,\Pi_{j = 1}^{N}(1-\beta_j\,z)}{z^N},
\end{equation*}
where $|\alpha_j|, |\beta_j| < e^{-A(\nu)}$ for all $j \in \{1,...,N\}$ and $c \not= 0$ is a constant. From this it is easy to see that 
$\varphi^{(N)}$ may be connected to $1$ through a continuous path in $GW_{\nu}$, i.e. $\varphi^{(N)} \in G_{0}W_{\nu}$. 
On the other hand, clearly 
\begin{equation*}
\lambda \varphi + (1-\lambda) \varphi^{(N)} = \varphi^{(N)} + \lambda (\varphi - \varphi^{(N)}), \quad \lambda \in [0,1],
\end{equation*}
connects $\varphi^{(N)}$ and $\varphi$ through a continuous path in $GW_{\nu}$ if $N$ is chosen sufficiently large, and so 
$\varphi \in G_{0}W_{\nu} = \exp W_{\nu}$.
We also mention the well-known fact that if $b$ belongs to a Banach algebra $\mathcal{B}$ and $f$ is a function analytic in a domain
containing $\sigma(b)$, then $f(b) \in \mathcal{B}$. 

Consequently, for $\varphi \in W_{\nu}$ with $\varphi \not= 0$  on $\mathcal{R}_{\nu}$, $wind(\varphi,0) = 0$, we have 
$\varphi_+, \varphi_-, \varphi_+^{-1}, \varphi_-^{-1} \in W_{\nu}$.    

We shall need some additional notation. Introduce, for $\varphi$ as above and $n \geq 0$, the semi-norms
\begin{equation}
  \label{Banachseminorms}
||\varphi||_{\nu,n} = \sum_{|k| \geq n} \nu_k\,|\varphi_k|, 
\end{equation}
and also write
\begin{equation}
|||\varphi|||_{\nu} = \max \left\{ ||\varphi_+||_{\nu}, ||\varphi_-||_{\nu}, ||\varphi_+^{-1}||_{\nu}, ||\varphi_-^{-1}||_{\nu}\right\}
\end{equation}
as well as
\begin{equation}
|||\varphi|||_{\nu,n} = \max \left\{ ||\varphi_+||_{\nu,n}, ||\varphi_-||_{\nu,n}, ||\varphi_+^{-1}||_{\nu,n}, ||\varphi_-^{-1}||_{\nu,n}\right\}.
\end{equation}
We will always replace $\nu$ by $0$ in \eqref{Banachnorm}, \eqref{Banachseminorms},... in case $\nu$ is the standard Wiener weight.
 
The first result in this paper is a new proof of the following basic theorem, which is essentially due to Widom. See \cite{BottSilb} for references and 
further discussion.
\begin{theorem}
  \label{mainresult}
Let $\nu$ be a Beurling weight. Suppose that $\varphi \in W_{\nu}$, that $\varphi(z) \ne 0$ for all $z \in \mathcal{R}_{\nu}$, and that 
$wind(\varphi,0) = 0$. Let $\varphi = \varphi_+\,\varphi_-$ be the Wiener-Hopf factorization of $\varphi$. 
Then $T_n(\varphi)$ is invertible for sufficiently large $n$, and there is a constant $c(\varphi)$ (independent of $n$) such that 
\begin{equation}
  \label{main1}
\left|T_n(\varphi)_{jk}^{-1} - T(\varphi)_{jk}^{-1}\right| \leq
c(\varphi) \cdot \min \left\{|||\varphi|||_{0,n+1-k}, |||\varphi|||_{0,n+1-j}\right\}
\end{equation}  
for $0 \leq j,k \leq n$. In particular, for any Beurling weight with $A(\nu) > 0$,
\begin{equation}
  \label{main2}
\left|T_n(\varphi)_{jk}^{-1} - T(\varphi)_{jk}^{-1}\right| \leq
c_{\nu}(\varphi) \cdot \min \left\{e^{-(n+1-k)\,A(\nu)}, e^{-(n+1-j)\,A(\nu)}\right\}.
\end{equation}
On the other hand, for Beurling weights which increase on $\mathbb{Z}_+$; $\nu_j \leq \nu_k\,$ for $0 \leq j < k$, 
\begin{equation}
  \label{main3}
\left|T_n(\varphi)_{jk}^{-1} - T(\varphi)_{jk}^{-1}\right| \leq
c_{\nu}(\varphi) \cdot \min \left\{\nu_{n+1-k}^{-1}, \nu_{n+1-j}^{-1}\right\}.
\end{equation}
\end{theorem}
\textit{Remarks.}
1. For symbols $\varphi$ which are positive on $\Gamma$ standard computations show that $T_n(\varphi)^{-1}$ exists for 
\textit{all} $n \geq 0$.\\
2. Of course, \eqref{main2} is true for all Beurling weights, but is only of interest if $A(\nu) > 0$.\\
This result has many applications. For a recent application to random growth models, see \cite{Joh}.

The second result concerns the relationship between the asymptotic properties of Verblunsky coefficients and the smoothness of the measures 
$d\mu$ on the unit circle. The result is the following extension of the $\mathcal{I}$-part of Baxter's theorem (see Section \ref{Baxtersection}). 
\begin{theorem}
  \label{Baxterextensiontheorem}
Let $\nu$ be a Beurling weight and $d\mu(z)= w(z)\,\frac{|dz|}{2\pi}$, a complex measure on the unit circle with the properties 
$w \in W_{\nu}$, $w(z) \not= 0$ for $z \in \mathcal{R}_{\nu}$ and $wind(w,0) = 0$. Then,
\begin{equation}
  \label{Baxterextension}
\sum_{n \geq n_0} \nu_n\,|\Phi_n(0)| < \infty,
\end{equation}
for some $n_0 = n_0(\nu)$ sufficiently large. 
\end{theorem}
As in the case of real weights, $\Phi_n = z^n + ...$ is the monic polynomial defined by the conditions
$\int_{\Gamma} \Phi_n(z)\,z^{-k}\,w(z)\,|dz| = 0$, $0 \leq k \leq n-1$. For complex-valued weights as above, such polynomials may not exist for all $n$.
However, for $n$ sufficiently large such polynomials exist and are unique. There are two ways to see this. Firstly, a simple computation
shows that polynomials $\Phi_n$ exist and are unique if the Toeplitz operator $\left(T_{n-1}(w)\right)_{0 \leq j,k \leq n-1}$ is invertible - but 
as remarked at the end of Section \ref{Errorsection} below this is true for $n$ sufficiently large. On the other hand, if the RHP in 
Section \ref{Baxtersection} below has a unique solution $Y$, then $Y_{11}$ is the desired (unique) polynomial. The existence of a unique solution 
$Y$ for $n$ sufficiently large is proven en route in the calculations of Section \ref{Baxtersection}. Of course, in case $w > 0$ 
(as in Baxter's theorem), the OPUC's $\Phi_n$ exist for all $n \geq 0$ and we take $n_0 = 0$ in \eqref{Baxterextension}.

Whereas the results (but not the methods!) mentioned above are basically classical, our third result, Theorem \ref{Simon-improvement}
given in Section \ref{Baxtersection}, is new. It is a further refinement of Baxter's theorem and may be regarded as a result about the 
Born approximation for a scattering problem on $\mathbb{Z}_+$. Together with results from Nevai and Totik \cite{NT}, one implication of this 
result is a strengthening (see Corollary \ref{Analyticity-corollary}) of an earlier result of Simon. As it turns out, Simon has now given an 
independent proof of this Corollary (see \cite{Sim3}).

Section \ref{IntegrableOp&RHPsection} briefly discusses techniques from the theory of integrable operators and RHP's which we will need in 
the sequel. Sections \ref{TruncatedToeplitz}-\ref{Errorsection} contain the proof of Theorem \ref{mainresult}. 
In Section \ref{Examplesection} we consider two examples illustrating the sharpness of the results in Section \ref{Baxtersection}. 
Finally, Section \ref{Ibragimovsection} contains a RH proof of the $\mathcal{I}$-part of a basic theorem of Golinskii-Ibragimov 
related to the Strong Szeg\"{o} Limit Theorem (see Theorem \ref{Ibragimovtheorem}, et seq.). For a proof of the Strong Szeg\"{o} Limit Theorem
based on RH techniques, we refer the reader to \cite{D1}.
 
\section{Integrable operators and Riemann-Hilbert problems.}
\label{IntegrableOp&RHPsection}
In this section we give a brief introduction to the theory of integrable operators and their connection to RHP's.
Let $\Sigma$ be an oriented contour in $\mathbb{C}$. We say that an operator $K$ acting in $L^2\left(\Sigma\right) = 
L^2\left(\Sigma, |dz|\right)$ is \textit{integrable} if it has a kernel of the form 
\begin{equation}
  \label{kernel}
K(z,z') = \frac{\sum_{j=1}^N f_j(z) g_j(z')}{z - z'}, \quad z, z' \in \Sigma,
\end{equation}
for some functions $f_i, g_j, 1 \leq i,j \leq N.$
The action of $K$ in $L^2 \left(\Sigma\right)$ is given by 
\begin{equation}
(Kh)(z) = i \pi \, \sum_{j=1}^N f_j(z) \left(H(hg_j)\right)(z), \quad h \in L^2\left(\Sigma\right), z \in \Sigma,
\end{equation}
where $H$ denotes the Hilbert-transform,
\begin{equation}
(Hh)(z) = \lim_{\epsilon \rightarrow 0}\,\frac{1}{i \pi} \int_{\left\{z' \in \Sigma\,:\,|z-z'| > \epsilon\right\}} \frac{h(z')}{z-z'} dz', 
\quad h \in L^2\left(\Sigma\right), z \in \Sigma. 
\end{equation}
In case the contour $\Sigma$ is such that the operator $H$ is bounded on $L^2 \left(\Sigma\right)$, and if $f_i, g_j \in 
L^{\infty}\left(\Sigma\right)$ for $1 \leq i,j \leq N$, then clearly $K$ defines a bounded operator on $L^2 \left(\Sigma\right)$.
Particular examples of integrable operators began to appear in the 1960's in the context of field theory and statistical models and  
some of the important elements of the general theory of such operators were present in the late 60's in \cite{Sak}, but the full theory 
of integrable operators as a distinguished class was presented only in the early 90's  in \cite{IIKS} (see also \cite{D1}).

Integrable operators have many remarkable properties, see \cite{IIKS}, \cite{D1}. 
In particular, if $K$ is an integrable operator with kernel as in \eqref{kernel}, with the property that $(1 - K)^{-1}$ exists,
and $(1 - K)^{-1} - 1 = R$ is also a kernel operator, then we learn from \cite{IIKS}, \cite{D1} that $R$ is also an integrable operator with kernel
\begin{equation}
  \label{Rkernel}
R(z,z') = \frac{\sum_{j=1}^N F_j(z) G_j(z')}{z - z'}, \quad z, z' \in \Sigma,
\end{equation}
where
\begin{equation}
F_i = (1 - K)^{-1}f_i, \quad G_i = (1 - K^T)^{-1}g_i, \quad 1\leq i \leq N.
\end{equation}
Moreover, (see \cite{IIKS}) these functions $F_i$ and $G_i$ can be computed in terms of a canonical auxiliary 
Riemann-Hilbert matrix factorization problem naturally associated with $K$, as described below.\\

We now recall the basic definition of a Riemann-Hilbert matrix factorization problem. Let $\Sigma$ be an oriented countour in $\mathbb{C}$,
as above. As we move along an arc in $\Sigma$ in the direction of the orientation we say, by convention, that the (+)-side (resp. (-)-side)  
lies to the left (resp. right). The data of a RHP consists of a pair $\left(\Sigma, v\right)$, where 
$v\,:\, \Sigma \rightarrow Gl\left(k,\mathbb{C}\right)$ and $v, v^{-1} \in L^{\infty}(\Sigma)$. In case $\Sigma$ is unbounded we demand 
that $v(z) \rightarrow I$ as $z \rightarrow \infty$. The (normalized) RHP consists in proving existence of a (unique) $k \times k$ matrix-function 
$m = m(z)$, known as the \textit{solution} of the RHP, satisfying\\
\quad $\bullet$ $m$ is analytic in $\mathbb{C}\backslash\Sigma$,\\
\quad $\bullet$ $m_+(z) = m_-(z)\,v(z)$, \quad $z \in \Sigma$,\\
\quad $\bullet$ $m(z) \rightarrow I$ as $z \rightarrow \infty$.\\
Here $m_\pm(z)$ denotes the limits of $m(z')$ as $z'$ approaches $z$ from the $(\pm)$-side of $\Sigma$. 
The matrix $v$ is called the \textit{jump matrix} for the RHP.
The precise sense in which the limits, $m_{\pm}(z) = \lim_{z' \rightarrow z} m(z')$ and $\lim_{z \rightarrow \infty} m(z) = I$, are attained 
is a technical matter (see e.g. \cite{CG} for details). The latter limit requires special care, in particular, when $\Sigma$ is unbounded. 
In all the RHP's that we consider in this paper, we will require in addition that\\
\quad $\bullet$ $m$ is continuous up to the boundary of $\mathbb{C} \backslash \Sigma$,\\
and also\\
\quad $\bullet$ $m(z) \rightarrow I$ uniformly as $z \rightarrow \infty$ in $\mathbb{C} \backslash \Sigma$.\\ 
The RHP $\left(\Sigma, v\right)$ reduces (see e.g. \cite{CG}) to the study of a singular integral operator on $\Sigma$ in the following way. 
Let
\begin{equation}
v(z) = \left(v_-(z)\right)^{-1}\left(v_+(z)\right), \quad z \in \Sigma,
\end{equation}
be any pointwise factorization of $v(z)$ with $v_\pm(z) \in Gl\left(k,\mathbb{C}\right)$. In case $\Sigma$ is unbounded we again
demand $v_\pm(z) \rightarrow I$ as $z \rightarrow \infty$. Define $\omega_\pm\,:\,\Sigma \rightarrow Gl\left(k,\mathbb{C}\right)$ through
the relations
\begin{equation}
v_\pm(z) = I \pm \omega_\pm(z), \quad z \in \Sigma.
\end{equation}
Denote the Cauchy operator by 
\begin{equation}
(Ch)(z) = \frac{1}{2 \pi i} \int_{\Sigma} \frac{h(z')}{z'-z} dz', \quad h \in L^2\left(\Sigma\right), z \in \mathbb{C} \backslash \Sigma,
\end{equation}
and set
\begin{equation}
\left(C_\pm h\right)(z) = \lim_{
\begin{array}{rl}
& \quad \quad z' \rightarrow z\\
&z'\in (\pm)\mbox{-side of } \Sigma
\end{array}
} (Ch)(z'), \quad h \in L^2\left(\Sigma\right), z \in \Sigma. 
\end{equation}
Standard computations show that 
\begin{equation}
C_\pm = \pm \frac{1}{2} - \frac{1}{2} H,
\end{equation}
so that 
\begin{equation}
  \label{cauchyidentities}
C_+ - C_- = 1, \quad C_+ + C_- = -H.
\end{equation}
For a given factorization $v = \left(I - \omega_-\right)^{-1} \left(I + \omega_+\right)$, define the operator
\begin{equation}
  \label{Comega}
C_{\omega} h = C_+\left(h\omega_-\right) + C_-\left(h\omega_+\right), 
\end{equation} 
for $k \times k$ matrix-valued functions $h$ in $L^2\left(\Sigma\right)$.  
Let $\mu \in I + L^2\left(\Sigma\right)$ be the solution of the singular integral equation 
\begin{equation}
  \label{muequation}
\left(1 - C_\omega\right) \mu = I.
\end{equation}

\textit{Remark.} For later purposes note that if $\Sigma$ is bounded, then $I \in L^2(\Sigma)$, and hence $\mu \in L^2(\Sigma)$.\\
Set 
\begin{equation}
  \label{generaltheory}
m(z) = I + C\left(\mu\left(\omega_+ + \omega_-\right)\right)(z), \quad z \in \mathbb{C} \backslash \Sigma.
\end{equation}
A basic computation using \eqref{cauchyidentities} and \eqref{muequation}, then shows that 
\begin{equation}
m_\pm(z) = \mu\,v_\pm, \quad z \in \Sigma.
\end{equation} 
Therefore, $m_+ = m_-\,v_-^{-1}\,v_+ = m_-\,v $. Clearly, $m$ is analytic in $\mathbb{C} \backslash \Sigma$ and $m(z) \rightarrow I$ as 
$z \rightarrow \infty$, so that, modulo technicalities, $m$ solves the RHP. Conversely, one verifies that if $m$ solves the RHP, 
then $\mu = m_+\,v_+^{-1} = m_-\,v_-^{-1}$ solves \eqref{muequation}. Thus, the existence (and uniqueness) of the solution of the RHP is 
equivalent to the existence (and uniqueness) of a solution $\mu \in I + L^2(\Sigma)$ of the singular 
integral equation \eqref{muequation} for any (and hence all) pointwise factorization(s) $v = \left(I - \omega_-\right)^{-1} \left(I + \omega_+\right)$.\\

We now return to our discussion of integrable operators. Suppose $K$ is an integrable operator with kernel as in \eqref{kernel}, and 
that $(1-K)^{-1}$ exists with $(1 - K)^{-1} - 1 = R$ also a kernel operator. The remarkable fact proven in \cite{IIKS}, \cite{D1} is the following:
the functions $F_i, G_i$ in the kernel \eqref{Rkernel} of the operator $R$ can be computed as
\begin{alignat}{4}
\label{Fform}
F &= (F_1,...,F_N)^T = \left(1 \mp i\pi\,f^T\,g\right)^{-1} m_\pm f,\\
\label{Gform}
G &= (G_1,...,G_N)^T = \left(1 \pm i\pi\,f^T\,g\right)^{-1} \big(m^T\big)_\pm^{-1} g,
\end{alignat}
where $m$ is the solution of the RHP $\left(\Sigma,v\right)$ with
\begin{equation}
  \label{vform}
v = I - \left(\frac{2 \pi i}{1+i\pi\,f^T\,g}\right)fg^T.
\end{equation}

\section{Truncated Toeplitz operators as integrable operators.}
\label{TruncatedToeplitz} 
From now on we will assume $\Gamma = \{z \in \mathbb{C} : |z| = 1\}$ to be oriented counterclockwise.
A direct calculation shows that for any polynomial $p = \sum_{j=0}^n a_j\,z^j \in \mathcal{P}_n$,
\begin{equation}
  \label{truncatedasintegrable}
\left(T_n p\right)(z) = \left(\left(1 - K_n\right)p\right)(z) = p(z) - \int_{\Gamma} K_n(z,z')p(z') dz',
\end{equation}
where $K_n = K_n(\varphi)\,:\,L^2\left(\Gamma\right) \rightarrow L^2\left(\Gamma\right)$ is the operator with kernel
\begin{equation}
  \label{kernelspec}
K_n(z,z') = \frac{z^{n+1}(z')^{-(n+1)}-1}{z-z'}\,\frac{1-\varphi(z')}{2\pi i}.
\end{equation}
Clearly, $K_n$ is an integrable operator on $L^2\left(\Gamma\right)$ of form \eqref{kernel}, where
\begin{alignat}{4}
f &= \left(f_1,f_2\right)^T = \left(z^{n+1},1\right)^T,\\
g &= \left(g_1,g_2\right)^T = \left(z^{-(n+1)}\frac{1-\varphi(z)}{2\pi i},-\frac{1-\varphi(z)}{2\pi i}\right)^T.
\end{alignat}
Since $f^T\,g\,= 0$ the formulas \eqref{Fform}, \eqref{Gform} and \eqref{vform} for the functions $F_i, G_j$ appearing in the kernel \eqref{Rkernel} 
of $R_n = \left(1-K_n\right)^{-1}-1$ simplify to 
\begin{equation}
  \label{FGform}
F = m_+\,f, \quad G = \left(m_+^T\right)^{-1}\,g,
\end{equation}
where $m$ solves the RHP $(\Gamma, v)$ with
\begin{equation}
  \label{vformspec}
v =
\begin{pmatrix}
\varphi & -z^{n+1}(\varphi-1)\\
z^{-(n+1)}(\varphi-1) & 2-\varphi
\end{pmatrix}.
\end{equation}
Clearly, 
\begin{equation*}
T_n(\varphi)\,z^l = \sum_{j=0}^{n} \varphi_{j-l}\,z^j, \quad 0 \leq l \leq n, 
\end{equation*}
and so (whenever $T_n$ is invertible) identity \eqref{truncatedasintegrable} implies: for $0 \leq l,k \leq n$
\begin{equation*}
2\pi\,\delta_{l,k} = \left(z^l, z^k\right)_{L^2\left(\Gamma,\,|dz|\right)} = \sum_{j=0}^{n} \varphi_{j-l} 
\left((1-K_n)^{-1}z^j, z^k\right)_{L^2\left(\Gamma,\,|dz|\right)}.
\end{equation*}
Hence,
\begin{equation}
  \label{outset1}
\left(T_n(\varphi)\right)_{j,k}^{-1} = \delta_{j,k} + \frac{1}{2\pi}\left(R_n(\varphi)\,z^k, z^j\right)_{L^2\left(\Gamma,\,|dz|\right)}, 
\quad 0 \leq j,k \leq n.
\end{equation}
This identity is basic for our proof of Theorem \ref{mainresult}. The invertibility of $T_n$, for large $n$, will be discussed below 
(see the end of Section \ref{Errorsection}).\\

In order to make the forthcoming ideas transparent, let us first assume that $\varphi$ is analytic in some annular domain 
$\left\{\rho < |z| < \rho^{-1}\right\}$, $0 < \rho < 1$. The basic observation is that the lower/upper factorization of $v$, which always exists:
\begin{equation}
  \label{vfactorization}
v = 
\begin{pmatrix}
1 & 0\\
z^{-(n+1)}(1-\varphi^{-1}) & 1
\end{pmatrix}
\begin{pmatrix}
\varphi & 0\\
0 & \varphi^{-1}
\end{pmatrix}
\begin{pmatrix}
1 & -z^{n+1}(1-\varphi^{-1})\\
0 & 1
\end{pmatrix},
\end{equation}
can then be analytically extended to the annulus.\\
Let $\rho < \rho^{(1)} < 1$. Define the function $m^{(1)}$ by
\begin{alignat}{4}
m^{(1)}(z) &= m(z), \quad |z| < \rho^{(1)},\\
m^{(1)}(z) &= m(z) \begin{pmatrix} 1 & -z^{n+1}(1-\varphi^{-1})\\0 & 1\end{pmatrix}^{-1}, \quad \rho^{(1)} < |z| < 1,\\
m^{(1)}(z) &= m(z) \begin{pmatrix} 1 & 0\\z^{-(n+1)}(1-\varphi^{-1}) & 1\end{pmatrix}, \quad 1 < |z| < (\rho^{(1)})^{-1},\\
m^{(1)}(z) &= m(z), \quad |z| > (\rho^{(1)})^{-1}.
\end{alignat} 
Then $m^{(1)}$ solves the RHP $\left(\Gamma^{(1)}, v^{(1)}\right)$, where 
$\Gamma^{(1)} = \left\{|z|=\rho^{(1)}\right\}\cup\left\{|z|=1\right\}\cup\left\{|z|=\left(\rho^{(1)}\right)^{-1}\right\}$, oriented
counterclockwise on each circle, and
\begin{alignat}{4}
v^{(1)}(z) &= \begin{pmatrix} 1 & -z^{n+1}(1-\varphi^{-1})\\0 & 1\end{pmatrix}, \quad |z| = \rho^{(1)},\\
v^{(1)}(z) &= \begin{pmatrix} \varphi & 0\\0 & \varphi^{-1}\end{pmatrix}, \quad |z| = 1,\\
v^{(1)}(z) &= \begin{pmatrix} 1 & 0\\z^{-(n+1)}(1-\varphi^{-1}) & 1\end{pmatrix}, \quad |z| = (\rho^{(1)})^{-1}.
\end{alignat}
As $n$ gets large, the solution $m^{(1)}$ of the RHP $\left(\Gamma^{(1)}, v^{(1)}\right)$ should (in some sense) be close to the solution
$m_{\infty}^{(1)}$ of the RHP $\left(\Gamma^{(1)}, v_{\infty}^{(1)}\right)$, where
\begin{alignat}{4}
\label{v1infty1}
v_{\infty}^{(1)}(z) &= I, \quad |z| = \rho^{(1)},\\
\label{v1infty2}
v_{\infty}^{(1)}(z) &= \begin{pmatrix} \varphi & 0\\0 & \varphi^{-1}\end{pmatrix}, \quad |z| = 1,\\
\label{v1infty3}
v_{\infty}^{(1)}(z) &= I, \quad |z| = (\rho^{(1)})^{-1}.
\end{alignat}
Standard computations show that the solution of \eqref{v1infty1}-\eqref{v1infty3} is given by
\begin{equation}
m_{\infty}^{(1)} = \exp\,\{C(\log(\varphi))\}^{\sigma_3},
\end{equation} 
where $\sigma_3 = \begin{pmatrix} 1 & 0\\0 & -1\end{pmatrix}$ denotes the third Pauli matrix.\\
Hence we expect that $m$ is close (in some sense) to $m_{\infty}$, where
\begin{alignat}{4}
m_{\infty}(z) &= m_{\infty}^{(1)}(z), \quad |z| < \rho^{(1)},\\
\label{interior}
m_{\infty}(z) &= m_{\infty}^{(1)}(z) \begin{pmatrix} 1 & -z^{n+1}(1-\varphi^{-1})\\0 & 1\end{pmatrix}, \quad \rho^{(1)} < |z| < 1,\\
m_{\infty}(z) &= m_{\infty}^{(1)}(z) \begin{pmatrix} 1 & 0\\z^{-(n+1)}(1-\varphi^{-1}) & 1\end{pmatrix}^{-1}, 
\quad 1 < |z| < (\rho^{(1)})^{-1},\\
m_{\infty}(z) &= m_{\infty}^{(1)}(z), \quad |z| > (\rho^{(1)})^{-1}.
\end{alignat} 
Finally, let us define 
\begin{equation}
  \label{Rinfinitykernel}
R_n^{\infty}(\varphi; z,z') = \frac{\sum_{j=1}^2 F_j^{\infty}(z) G_j^{\infty}(z')}{z - z'}, \quad z, z' \in \Gamma,
\end{equation}
where
\begin{equation}
  \label{FinfGinfform}
F^{\infty}(z) = m_{\infty,+}(z)\,f(z), \quad G^{\infty}(z) = \left(m_{\infty,+}^T\right)^{-1}(z)\,g(z), \quad z \in \Gamma,
\end{equation}
and also write
\begin{equation}
  \label{outset2}
\left(T_n^{\infty}(\varphi)\right)_{j,k}^{-1} = \delta_{j,k} + \frac{1}{2\pi}\left(R_n^{\infty}(\varphi)\,z^k, z^j\right)_{L^2\left(\Gamma,\,|dz|\right)}.
\end{equation}
We emphasize that we use the left-hand side of \eqref{outset2} only as a formal symbol for the quantity on the right-hand side. 
By the above consideration, we expect
\begin{equation}
  \label{expectation}
\left(T_n(\varphi)\right)_{jk}^{-1} \sim \left(T_n^{\infty}(\varphi)\right)_{jk}^{-1}.
\end{equation}

Although in this section we have assumed analyticity of $\varphi$ in order to motivate our calculations, note the following: even in case 
that $\varphi$ is \textit{not} analytic in an annulus we still define 
$m_{\infty,+}(z) = m_{\infty,+}^{(1)}(z) \begin{pmatrix} 1 & -z^{n+1}(1-\varphi^{-1})\\0 & 1\end{pmatrix}$, $z \in \Gamma$, and also  
$F^{\infty}, G^{\infty}, R_n^{\infty}$ and $\left(T_n^{\infty}\right)_{j,k}^{-1}$ in the same way.   
Under the only assumption that $\varphi$ belongs to $W_{\nu}$ we still expect \eqref{expectation} to be true.

\textit{Remark.} In case $\varphi$ is analytic in an annulus, $m_{\infty,+}$ is the boundary value on $\Gamma$ of a piecewise analytic
function $m_{\infty}$ which solves a RHP. In general, for $\varphi \in W_{\nu}$, this is no longer true. 
\newpage
\section{Explicit computation of $\left(T_n^{\infty}(\varphi)\right)_{j,k}^{-1}$ for $\varphi$ in $W_{\nu}$.}
Solving \eqref{v1infty1}-\eqref{v1infty3} for $m_{\infty}^{(1)}$ and using definition \eqref{Rinfinitykernel} and 
the Wiener-Hopf factorization $\varphi = \varphi_+\,\varphi_-$, we obtain 
\begin{alignat}{4}
  \label{Rninfty}
R_n^{\infty}(\varphi; z,z') &= \frac{1}{2\pi i} \frac{1}{(z')^{n+1}}\frac{1}{z-z'}
\Bigg[(z')^{n+1}\frac{\varphi_+(z')}{\varphi_+(z)} - z^{n+1}\frac{\varphi_-(z')}{\varphi_-(z)}\\
\notag
&+z^{n+1}\frac{1}{\varphi_-(z)\varphi_+(z')} - (z')^{n+1}\frac{1}{\varphi_+(z)\varphi_-(z')}
\Bigg]. 
\end{alignat}
In order to evaluate the right-hand side of \eqref{outset2} further, it is convenient to assume again that $\varphi$ is analytic in an 
annulus $\left\{\rho < |z| < \rho^{-1}\right\}$, $0 < \rho < 1$. Clearly then $\varphi_{\pm}$ are also analytic in the same
annulus. We will later remove this analyticity assumption (see below).
Writing $\Gamma_{\epsilon} = \left\{z \in \mathbb{C}\,:\,|z| = 1 - \epsilon\right\}$, $\epsilon > 0$ sufficiently small,  
and using Cauchy's theorem as well as the elementary identity 
\begin{equation*}
\frac{1}{z-z'} = \frac{1}{z} \sum_{m=0}^{\infty} \left(\frac{z'}{z}\right)^m, \quad z \in \Gamma,\,z' \in \Gamma_{\epsilon},
\end{equation*} 
we then obtain for $0 \leq j,k \leq n$ 
\begin{alignat*}{4}
\frac{1}{2\pi} &\left(R_n^{\infty}(\varphi)\,z^k, z^j\right)_{L^2\left(\Gamma,\,|dz|\right)} = 
\lim_{\epsilon \downarrow 0}\,\frac{1}{2\pi} \int_{\Gamma}\left(\int_{\Gamma_{\epsilon}} R_n^{\infty}(z,z') (z')^k\,dz'\right) z^{-j}\,\frac{dz}{iz}\\
&= \frac{1}{(2\pi i)^2} \lim_{\epsilon \downarrow 0}\,\sum_{m=0}^{\infty}
\Bigg[
\int_{\Gamma_{\epsilon}}(z')^{k+m} \varphi_+(z')dz' \cdot \int_{\Gamma} z^{-(j+2+m)} \frac{1}{\varphi_+(z)}dz\\
&-\int_{\Gamma_{\epsilon}}(z')^{k-n-1+m} \varphi_-(z')dz' \cdot \int_{\Gamma} z^{-(j+2+m-n-1)} \frac{1}{\varphi_-(z)}dz\\  
&+\int_{\Gamma_{\epsilon}}(z')^{k-n-1+m} \frac{1}{\varphi_+(z')}dz' \cdot \int_{\Gamma} z^{-(j+2+m-n-1)} \frac{1}{\varphi_-(z)}dz\\
&-\int_{\Gamma_{\epsilon}}(z')^{k+m} \frac{1}{\varphi_-(z')}dz' \cdot \int_{\Gamma} z^{-(j+2+m)} \frac{1}{\varphi_+(z)}dz
\Bigg]\\
&=\sum_{m=0}^{\infty} 
\Bigg[
\left(\varphi_+\right)_{-k-1-m} \left(\varphi_+^{-1}\right)_{j+1+m}
-\left(\varphi_-\right)_{n-k-m} \left(\varphi_-^{-1}\right)_{j-n+m}\\
&+\left(\varphi_+^{-1}\right)_{n-k-m} \left(\varphi_-^{-1}\right)_{j-n+m}  
-\left(\varphi_-^{-1}\right)_{-k-1-m} \left(\varphi_+^{-1}\right)_{j+1+m} 
\Bigg]\\
&=0 - \left[T\left(\varphi_-\right) T\left(\varphi_-^{-1}\right)\right]_{n-k,n-j} +
\sum_{m=j+k-n}^{\infty} \left(\varphi_+^{-1}\right)_{j-m} \left(\varphi_-^{-1}\right)_{m-k}\\
&-\sum_{m=-\infty}^{-1} \left(\varphi_+^{-1}\right)_{j-m} \left(\varphi_-^{-1}\right)_{m-k}.
\end{alignat*}
Since $\left[T\left(\varphi_-\right) T\left(\varphi_-^{-1}\right)\right]_{n-k,n-j} = \delta_{j,k}$ we obtain upon insertion
into \eqref{outset2}: 
\begin{alignat}{4}
\notag
\left(T_n^{\infty}(\varphi)\right)_{j,k}^{-1} &= \sum_{m=j+k-n}^{\infty} \left(\varphi_+^{-1}\right)_{j-m} \left(\varphi_-^{-1}\right)_{m-k}
-\sum_{m=-\infty}^{-1} \left(\varphi_+^{-1}\right)_{j-m} \left(\varphi_-^{-1}\right)_{m-k}\\
\label{explicitformula}
&=\left[T\left(\varphi_+^{-1}\right) T\left(\varphi_-^{-1}\right)\right]_{j,k} - 
\sum_{m=n+1-j-k}^{\infty} \left(\varphi_+^{-1}\right)_{j+m} \left(\varphi_-^{-1}\right)_{-(m+k)}. 
\end{alignat}
As we shall now see, the basic identity \eqref{explicitformula} remains valid if we only assume that $\varphi \in W_{\nu}$, i.e. without the 
restriction that $\varphi$ be analytic in an annular neighborhood of the unit circle.
To see this, let us write $\varphi = e^w$, where $w(z) = \sum_{-\infty}^{\infty} w_j\,z^j$, $z \in \Gamma$. Put 
$\varphi^{(N)} = e^{w^{(N)}}$, where $w^{(N)}(z) =  \sum_{-N}^{N} w_j\,z^j$, $z \in \Gamma$. Then $w, w^{(N)} \in W_{\nu}$. 
Observe that 
\begin{equation}
  \label{Linftyconvergence}
\lim_{N \rightarrow \infty}\,\left|\left|\left(\varphi_{\pm}\right)^{\pm 1} - \left(\varphi_{\pm}^{(N)}\right)^{\pm 1}
\right|\right|_{L^{\infty}\left(\Gamma\right)} = 0.   
\end{equation}
For instance, writing $w_{+}(z) = \sum_{0}^{\infty} w_j\,z^j$ and $w_{+}^{(N)}(z) = \sum_{0}^{N} w_j\,z^j$, we have  
\begin{equation*}
\varphi_{+} - \varphi_{+}^{(N)} = e^{w_{+}}-e^{w_{+}^{(N)}} = e^{w_{+}^{(N)}}\,\cdot\,\left(e^{\widetilde{w}_N} - 1\right),
\end{equation*}
where $\widetilde{w}_N(z) = \sum_{N+1}^{\infty} w_j\,z^j$. 
On the other hand,
\begin{equation*}
\left|e^{\widetilde{w}_N}-1\right| = \left|\int_0^1 \frac{d}{dt} e^{t\,\widetilde{w}_N}\,dt\right| = 
\left|\int_0^1 \widetilde{w}_N\,e^{t\,\widetilde{w}_N}\,dt\right| \leq \left|\widetilde{w}_N\right|\,\max_{0 \leq t \leq 1}\, 
\left|e^{t\,\widetilde{w}_N}\right|,
\end{equation*}
and since 
\begin{equation*}
\left|\widetilde{w}_N(z)\right| = \left|\sum_{N+1}^{\infty} w_j\,z^j\right| \leq \sum_{N+1}^{\infty} \left|w_j\right|, \quad z \in \Gamma,
\end{equation*}
the statement \eqref{Linftyconvergence} clearly follows in this case from the fact that $w \in W_{\nu}$. The other cases are almost 
identical. Since $\varphi^{(N)}$ is obviously analytic in $\mathbb{C} \backslash \{0\}$ the identity \eqref{explicitformula} is valid with 
$\varphi$ replaced by $\varphi^{(N)}$.  We shall now see that each term converges as $N \rightarrow \infty$
to the same term with $\varphi$. Firstly, 
\begin{equation}
  \label{term1}
\lim_{N \rightarrow \infty}\,\left(T_n^{\infty}(\varphi^{(N)})\right)_{j,k}^{-1} = \left(T_n^{\infty}(\varphi)\right)_{j,k}^{-1}.
\end{equation}
To see why, note from formula \eqref{Rninfty} that the operator 
$R_n^{\infty}(\varphi)$ consists of four parts, all being of the form $\psi_j\,H\,\chi_j$, $j = 1,..,4$. Here $\psi_i$, $\chi_j\,:\,
L^2(\Gamma) \rightarrow L^2(\Gamma)$, $1 \leq i,j \leq 4$,  are operators of multiplication. For instance (ignoring a factor 2), 
$\psi_1$ is multiplication by $\varphi_+^{-1}$ and $\chi_1$ is multiplication by $\varphi_+$. Using 
\eqref{Linftyconvergence} and $L^2$-boundedness of $H$ one therefore sees that 
\begin{equation*}
\lim_{N \rightarrow \infty}\,\left|\left|R_n^{\infty}(\varphi)-R_n^{\infty}\left(\varphi^{(N)}\right)\right|\right|_{L^2 \rightarrow L^2} = 0,
\end{equation*}
so that \eqref{term1} follows from \eqref{outset2}. Secondly, that
\begin{equation*}
\lim_{N \rightarrow \infty}\,\left[T\left((\varphi_+^{(N)})^{-1}\right) T\left((\varphi_-^{(N)})^{-1}\right)\right]_{j,k} = 
\left[T\left(\varphi_+^{-1}\right) T\left(\varphi_-^{-1}\right)\right]_{j,k},
\end{equation*}
follows similarly from \eqref{Linftyconvergence} and the basic estimates
\begin{equation*}
\left|\left|T\left(\varphi_{\pm}^{-1}\right) - T\left((\varphi_{\pm}^{(N)})^{-1}\right)\right|\right|_{L^2 \rightarrow L^2} \leq 
\left|\left|\left(\varphi_{\pm}\right)^{-1} - (\varphi_{\pm}^{(N)})^{-1}\right|\right|_{L^{\infty}\left(\Gamma\right)}.
\end{equation*}
Finally, we have 
\begin{alignat}{4}
  \label{term3}
\lim_{N \rightarrow \infty}\,&\bigg\{\sum_{m=n+1-j-k}^{\infty} \left((\varphi_+^{(N)})^{-1}\right)_{j+m}
\left((\varphi_-^{(N)})^{-1}\right)_{-(m+k)} \bigg\} \\
\notag
&= \sum_{m=n+1-j-k}^{\infty} \left(\varphi_+^{-1}\right)_{j+m} \left(\varphi_-^{-1}\right)_{-(m+k)}. 
\end{alignat}
To see why, first note that by a computation almost identical to that giving the inequality \eqref{errorestimate} below, we immediately
obtain 
\begin{alignat*}{4}
&\Bigg|\sum_{m=n+1-j-k}^{\infty} \left(\varphi_+^{-1}\right)_{j+m} \left(\varphi_-^{-1}\right)_{-(m+k)}\\
\notag
&-\sum_{m=n+1-j-k}^{\infty} \left((\varphi_+^{(N)})^{-1}\right)_{j+m}
\left((\varphi_-^{(N)})^{-1}\right)_{-(m+k)}\Bigg|\\
&\leq \bigg(\left|\left|\varphi_+^{-1}-(\varphi_+^{(N)})^{-1}\right|\right|_{\nu}\,
\left|\left|\varphi_-^{-1}\right|\right|_{\nu} + \left|\left|(\varphi_+^{(N)})^{-1}\right|\right|_{\nu}\,
\left|\left|\varphi_-^{-1}-(\varphi_-^{(N)})^{-1}\right|\right|_{\nu}
\bigg).
\end{alignat*}
On the other hand, with $w_{-}(z) = \sum_{-\infty}^{-1} w_j\,z^j$ and $w_{-}^{(N)}(z) = \sum_{-N}^{-1} w_j\,z^j$, we get
\begin{alignat*}{4}
&\left|\left|\varphi_{\pm}^{-1}-(\varphi_{\pm}^{(N)})^{-1}\right|\right|_{\nu} = 
\left|\left|e^{-w_{\pm}}-e^{-w_{\pm}^{(N)}}\right|\right|_{\nu}\\
&=\left|\left|\sum_{k \in \mathbb{N}} \frac{(-1)^k}{k!}\,
\left((w_{\pm})^k-(w_{\pm}^{(N)})^k\right)\right|\right|_{\nu}\\
&= \left|\left|\sum_{k \in \mathbb{N}} \frac{(-1)^k}{k!}\,\left(w_{\pm}-w_{\pm}^{(N)}\right)\,\sum_{j=0}^{k-1}\, 
(w_{\pm}^{(N)})^{k-1-j}\,\left(w_{\pm}\right)^j\right|\right|_{\nu}\\
&\leq \left|\left|w_{\pm}-w_{\pm}^{(N)}\right|\right|_{\nu}\,\sum_{k \in \mathbb{N}} \frac{1}{k!}\,k\,
\left|\left|w_{\pm}\right|\right|_{\nu}^{k-1} = \left|\left|w_{\pm}-w_{\pm}^{(N)}\right|\right|_{\nu}\,
\exp\left(||w_{\pm}||_{\nu}\right),
\end{alignat*}
since $||\cdot||_{\nu}$ is submultiplicative and $\left|\left|w_{\pm}^{(N)}\right|\right|_{\nu} \leq 
\left|\left|w_{\pm}\right|\right|_{\nu}$.
Obviously $w_{\pm}^{(N)} \rightarrow w_{\pm}$ in $W_{\nu}$ as $N \rightarrow \infty$, which completes the proof of \eqref{term3}.\\ 
From now on all assumptions of analyticity will be dropped, and from this section we shall only keep the basic fact that identity 
\eqref{explicitformula} is valid for all $\varphi \in W_{\nu}$.
\newpage
\section{Estimates of the remainder.}
\label{Errorsection}
In this section we shall provide the necessary estimates of the remainder. Assume that $\varphi \in W_{\nu}$, $\varphi \neq 0$ 
on $\mathcal{R}_{\nu}$ and $wind(\varphi,0)=0$. Then $\varphi_{\pm} \in W_{\nu}$ and $\varphi_{\pm}^{-1} \in W_{\nu}$. 
Clearly, for $0 \leq j,k \leq n$, we have
\begin{equation*}
\left|\left(\varphi_+^{-1}\right)_{j+m}\right| \leq \sum_{l=n+1-k}^{\infty} \left|\left(\varphi_+^{-1}\right)_l\right|,\quad m \geq n+1-j-k. 
\end{equation*} 
We can therefore estimate the ``error term'' in \eqref{explicitformula} as follows; for $0 \leq j, k \leq n$ 
\begin{alignat}{4}
  \label{errorestimate}
\Bigg|&\sum_{m=n+1-j-k}^{\infty} \left(\varphi_+^{-1}\right)_{j+m} \left(\varphi_-^{-1}\right)_{-(m+k)}\Bigg| \leq
\\
\notag
&\leq \left(\sum_{l=n+1-k}^{\infty} \left|\left(\varphi_+^{-1}\right)_l\right|\right)\cdot
\left(\sum_{l=n+1-j}^{\infty} \left|\left(\varphi_-^{-1}\right)_{-l}\right|\right)\\
\notag
&\leq |||\varphi|||_0 \cdot \min \left\{|||\varphi|||_{0,n+1-k}, |||\varphi|||_{0,n+1-j}\right\}. 
\end{alignat}
The main part of the proof of Theorem \ref{mainresult}, namely that of inequality \eqref{main1}, is complete once we prove that the estimate:
\begin{equation}
\left|\left(R_n\,z^k, z^j\right)_{L^2\left(\Gamma\right)}-\left(R_n^{\infty}\,z^k, z^j\right)_{L^2\left(\Gamma\right)}\right| \leq 
c(\varphi) \cdot |||\varphi|||_{0,n+1} 
\end{equation}
is valid for $0 \leq j,k \leq n$, with $c(\varphi)$ independent of $n$ (for $n$ sufficiently large).\\
First note that (see \eqref{FGform}) 
\begin{alignat*}{4}
\left(R_n\,z^k, z^j\right)_{L^2\left(\Gamma\right)} &= 
\int\int_{\Gamma \times \Gamma} \frac{F^T(z) G(z')}{z-z'}\,(z')^k\,z^{-j}\,dz'\,\frac{dz}{iz}\\
&=\lim_{\epsilon \downarrow 0} \int\int_{\left\{(z,z') \in \Gamma \times \Gamma\,:\,|z-z'| > \epsilon\right\}} 
\frac{F^T(z) G(z')}{z-z'}\,(z')^k\,z^{-j}\,dz'\,\frac{dz}{iz}\\
&=\pi\,\int_{\Gamma}\,F^{T}(z) z^{-(j+1)}\left(\lim_{\epsilon \downarrow 0}\int_{|z-z'| > \epsilon}\frac{G(z')(z')^k}{z-z'}\frac{dz'}{i\pi}\right)dz\\
&=\pi\,\int_{\Gamma}\,F^{T}(z) z^{-(j+1)}\,H\left(G(\diamond)\,\diamond^k\right)(z)dz\\ 
&=i\pi\,\int_{\Gamma} H\left(G^T(\diamond)\,\diamond^k\right)(z)\,\overline{\overline{F(z)}\,z^j}\,\frac{dz}{iz}\\ 
&= i\pi\,\left(H\left(G(\diamond)\,\diamond^k\right), \overline{F(\diamond)}\,\diamond^j\right)_{L^2\left(\Gamma\right)}.
\end{alignat*}
Therefore (see \eqref{FinfGinfform}), 
\begin{alignat*}{4}
&\left|\left(R_n\,z^k, z^j\right)_{L^2\left(\Gamma\right)} - \left(R_n^{\infty}\,z^k, z^j\right)_{L^2\left(\Gamma\right)} \right| =\\
&=\left|i\pi\,\left(H\left(G(\diamond)\,\diamond^k\right), \overline{F(\diamond)}\,\diamond^j\right)_{L^2\left(\Gamma\right)}-
i\pi\,\left(H\left(G^{\infty}(\diamond)\,\diamond^k\right), \overline{F^{\infty}(\diamond)}\,\diamond^j\right)_{L^2\left(\Gamma\right)}\right|\\
&\leq c\,\cdot\,\left(||G-G^{\infty}||_{L^2\left(\Gamma\right)}\,\cdot\,||F||_{L^2\left(\Gamma\right)} 
+ ||G^{\infty}||_{L^2\left(\Gamma\right)}\,\cdot\,||F-F^{\infty}||_{L^2\left(\Gamma\right)}\right),
\end{alignat*}
by $L^2$-boundedness of $H$.\\
We shall now prove that 
\begin{alignat}{4}
\label{keyerrorestimate1}
||G-G^{\infty}||_{L^2\left(\Gamma\right)}\,,\,||F-F^{\infty}||_{L^2\left(\Gamma\right)} &\leq c(\varphi) \cdot |||\varphi|||_{0,n+1},\\
\label{keyerrorestimate2}
||F||_{L^2\left(\Gamma\right)}\,,\,||G^{\infty}||_{L^2\left(\Gamma\right)} &\leq c(\varphi),
\end{alignat}
with $c(\varphi)$ independent of $n$ (for $n$ sufficiently large). For this we need the following elementary lemma:
\begin{lemma}
  \label{decaylemma}
For $n \geq 0$ and $f \in W_{\nu}$, 
\begin{equation}
||C_+\left(z^{-n}\,f\right)||_{L^2\left(\Gamma\right)} \leq \sqrt{2\pi}\,\sum_{k=n}^{\infty} \left|f_k\right|, \quad 
||C_-\left(z^n\,f\right)||_{L^2\left(\Gamma\right)} \leq \sqrt{2\pi}\,\sum_{k=n+1}^{\infty} \left|f_{-k}\right|.
\end{equation}
\end{lemma}
\textit{Proof.} We shall prove only the first bound, since the other is almost identical.\\ 
It is easy to verify (and we have already used several times without notice) the fact that $C_+$ agrees with the Riesz projection 
$P_+ : L^2(\Gamma) \rightarrow H_+$ on $L^2(\Gamma)$.  
Thus, 
\begin{equation*}
C_+\left(z^{-n}\,f\right)(z) = \sum_{k=0}^{\infty}\widehat{z^{-n}\,f}\,(k)\,z^k\, = \sum_{k=0}^{\infty} f_{k+n}\,z^k, 
\end{equation*}
so by Parseval
\begin{equation*}
||C_+\left(z^{-n}\,f\right)||_{L^2\left(\Gamma\right)}^2 = \sqrt{2\pi}\,\sum_{k=n}^{\infty} \left|f_k\right|^2 \leq 
\sqrt{2\pi}\,\left(\sum_{k=n}^{\infty} \left|f_k\right|\right)^2.
\end{equation*}$\Box$\\
The estimates \eqref{keyerrorestimate1} follow from the inequality
\begin{equation}
  \label{keyerrorestimate3}
||m_+-m_{\infty,+}||_{L^2\left(\Gamma\right)} \leq c(\varphi) \cdot |||\varphi|||_{0,n+1},
\end{equation}
with $c(\varphi)$ independent of $n$ (for $n$ sufficiently large), which we shall now prove. In view of \eqref{vfactorization} it is natural to put 
\begin{equation*}
\delta = \exp \left\{C(\log \varphi)\right\}, \quad \delta_{\pm} = \exp \left\{C_{\pm}(\log \varphi)\right\},
\end{equation*}
and 
\begin{equation*}
M = m\,\delta^{-\sigma_3}, 
\end{equation*}
where again $\sigma_3 = \begin{pmatrix} 1 & 0\\0 & -1\end{pmatrix}$ denotes the third Pauli matrix. Note that $\delta_+ = \varphi_+$ and 
$\delta_- = \varphi_-^{-1}$. Then, 
\begin{equation*}
M_+ = M_-\,v^M, 
\end{equation*}
where $v^M = \delta_{-}^{\sigma_3}\,v\,\delta_{+}^{-\sigma_3}$. A computation gives that
\begin{equation*}
v^M = \left(v_-^M\right)^{-1}\,v_+^M = \left(I-\omega_-^M\right)^{-1}\,\left(I+\omega_+^M\right),
\end{equation*}
where
\begin{equation}
  \label{omegaMplusminus}
\omega_-^M =
\begin{pmatrix}
0 & 0\\
z^{-(n+1)}\left(1-\varphi^{-1}\right)\,\delta_-^{-2} & 0
\end{pmatrix}, \quad 
\omega_+^M =
\begin{pmatrix}
0 & -z^{n+1}\left(1-\varphi^{-1}\right)\,\delta_+^2\\
0 & 0
\end{pmatrix}.
\end{equation}
We know that 
\begin{equation*}
M_{\pm} = \mu^M\,v_{\pm}^M,
\end{equation*}
where
\begin{equation}
  \label{muMequation}
\left(1-C_{\omega^M}\right) \mu^M = I, \quad \mu^{M} \in L^2(\Gamma). 
\end{equation}
Hence, 
\begin{equation}
\label{m+frommuM}
m_{\pm} = M_{\pm}\,\delta_{\pm}^{\sigma_3} = \mu^M\,v_{\pm}^M\,\delta_{\pm}^{\sigma_3},
\end{equation}
with $\mu^M$ given as the solution of the singular integral equation \eqref{muMequation}.
Also,
\begin{alignat}{4}
\label{m+inftyfrommuM}
m_{\infty,+}(z) &= m_{\infty,+}^{(1)}\,(z)\,\begin{pmatrix}
1 & -z^{n+1}\left(1-\varphi^{-1}\right)\\
0 & 1
\end{pmatrix}
=\\
\notag
&=\delta_+^{\sigma_3}\,(z)\,\begin{pmatrix}
1 & -z^{n+1}\left(1-\varphi^{-1}\right)\\
0 & 1
\end{pmatrix}
= v_+^M\,\delta_+^{\sigma_3},
\end{alignat}
for $z \in \Gamma$. Combining \eqref{m+frommuM} and \eqref{m+inftyfrommuM} we see that 
\begin{equation}
  \label{mminusminftyidentity}
m_+ - m_{\infty,+} = \left(\mu^M-I\right)\,v_+^M\,\delta_+^{\sigma_3}.
\end{equation}
On the other hand,
\begin{equation}
\mu^M-I=\left(1-C_{\omega^M}\right)^{-1}\,C_{\omega^M}\,I=\left(1-C_{\omega^M}\right)^{-1}\left(C_+\,\omega_-^M+C_-\,\omega_+^M\right).
\end{equation}
By Lemma \ref{decaylemma}
\begin{alignat}{4}
\notag
\left|\left|C_+\,\omega_-^M\right|\right|_{L^2\left(\Gamma\right)} &= 
\left|\left|C_+\left(z^{-(n+1)}\,\varphi^{-1}\,\delta_-^{-2}\right)\right|\right|_{L^2\left(\Gamma\right)} \leq 
\sqrt{2\pi}\,\sum_{k=n+1}^{\infty} |(\delta_+^{-1}\,\delta_-^{-1})_k| \\
&\leq \sqrt{2\pi}\,\sum_{l=0}^{\infty} |(\delta_-^{-1})_{-l}| \cdot \sum_{l=n+1}^{\infty} |(\delta_+^{-1})_{l}| 
\leq \sqrt{2\pi}\,|||\varphi|||_0 \cdot |||\varphi|||_{0,n+1}.
\end{alignat}
Similarly,
\begin{equation}
\left|\left|C_-\,\omega_+^M\right|\right|_{L^2\left(\Gamma\right)} \leq \sqrt{2\pi}\,|||\varphi|||_0 \cdot |||\varphi|||_{0,n+2}.
\end{equation}
Furthermore, clearly
\begin{equation}
  \label{Linftyestimate}
\left|\left|v_+^M\,\delta_+^{\sigma_3}\right|\right|_{L^{\infty}\left(\Gamma\right)} \leq 4\,|||\varphi|||_0. 
\end{equation}
Combining \eqref{mminusminftyidentity}-\eqref{Linftyestimate} we see that the proof of inequality 
\eqref{keyerrorestimate3} is complete once we show that $\left(1-C_{\omega^M}\right)^{-1}$
exists for $n$ sufficiently large, and that 
\begin{equation}
  \label{Inversebound}
\left|\left|\left(1-C_{\omega^M}\right)^{-1}\right|\right|_{L^2\left(\Gamma\right) \rightarrow L^2\left(\Gamma\right)} \leq c(\varphi),
\end{equation}
for $n$ sufficiently large, with $c(\varphi)$ independent of $n$.
One sees that existence of $\left(1-C_{\omega^M}^2\right)^{-1}$ implies existence of $\left(1-C_{\omega^M}\right)^{-1}$ and that 
\begin{equation}
  \label{inverseidentity}
\left(1-C_{\omega^M}\right)^{-1} = \left(1+C_{\omega^M}\right)\,\left(1-C_{\omega^M}^2\right)^{-1},
\end{equation}
whenever both inverses exists. But it is not difficult to see that 
\begin{equation}
  \label{C2limit1}
\lim_{n \rightarrow \infty}\,\left|\left|C_{\omega^M}^2\right|\right|_{L^2\left(\Gamma\right) \rightarrow L^2\left(\Gamma\right)} = 0.
\end{equation}
To see this, introduce the abbreviations $a = \left(1-\varphi^{-1}\right)\,\delta_-^{-2}$, $\alpha_n(z) = a\,z^{-(n+1)}$, 
$b = -\left(1-\varphi^{-1}\right)\,\delta_+^{-2}$ and $\beta_n(z) = b\,z^{n+1}$. A direct computation gives that
\begin{equation*}
C_{\omega^M}^2\,h =  
\begin{pmatrix}
C_+\left(\alpha_n\,C_-\left(\beta_n\,h_{11}\right)\right) & C_-\left(\beta_n\,C_+\left(\alpha_n\,h_{12}\right)\right)\\
C_+\left(\alpha_n\,C_-\left(\beta_n\,h_{21}\right)\right) & C_-\left(\beta_n\,C_+\left(\alpha_n\,h_{22}\right)\right)
\end{pmatrix},
h =
\begin{pmatrix}
h_{11} & h_{12}\\
h_{21} & h_{22}
\end{pmatrix}. 
\end{equation*} 
Consider $C_+\left(\alpha_n\,C_-\left(\beta_n\,h_{11}\right)\right)$, say. Obviously,
\begin{equation*}
\left|\left|C_+\left(\alpha_n\,C_-\left(\beta_n\,h_{11}\right)\right)\right|\right|_{L^2\left(\Gamma\right)} \leq 
\left|\left|C_+\,\alpha_n\,C_-\right|\right|_{L^2\left(\Gamma\right) \rightarrow L^2\left(\Gamma\right)}\,
\left|\left|\beta_n\right|\right|_{L^{\infty}\left(\Gamma\right)}\,\left|\left|h_{11}\right|\right|_{L^2\left(\Gamma\right)}.
\end{equation*} 
But clearly 
\begin{equation}
  \label{C2limit2}
\lim_{n \rightarrow \infty}\,\left|\left|C_+\,\alpha_n\,C_-\right|\right|_{L^2\left(\Gamma\right) \rightarrow L^2\left(\Gamma\right)} = 0.
\end{equation}
To see this, let $\epsilon > 0$. For $N$ sufficiently large
\begin{equation*}
\bigg|\bigg|a - \sum_{|k| \leq N} a_k\,z^k\bigg|\bigg|_{L^{\infty}\left(\Gamma\right)} < \epsilon.
\end{equation*}
Put $\widetilde{a}(z) = \sum_{|k| \leq N} a_k\,z^k$.
Then,
\begin{alignat*}{4}
&\left|\left|C_+\,\alpha_n\,C_-\right|\right|_{L^2 \rightarrow L^2} \leq\\
&\leq \left|\left|C_+\,(a-\widetilde{a})\,z^{-(n+1)}\,C_-\right|\right|_{L^2 \rightarrow L^2}+
\left|\left|C_+\,\widetilde{a}\,z^{-(n+1)}\,C_-\right|\right|_{L^2 \rightarrow L^2} \leq\\
&\leq ||C_+||_{L^2 \rightarrow L^2}\,||a-\widetilde{a}||_{L^{\infty}}\,
||C_-||_{L^2 \rightarrow L^2} + \bigg|\bigg|C_+\sum_{|k| \leq N} a_k\,z^{-(n-k+1)}\,C_-\bigg|\bigg|_{L^2 \rightarrow L^2}.
\end{alignat*}
Clearly the first term is $\epsilon$-small, whereas the second is zero for $n > N-2$. This verifies \eqref{C2limit2} and therefore
\eqref{C2limit1}. 
Using \eqref{inverseidentity}, \eqref{C2limit1} we immediately obtain
\begin{equation*}
\left|\left|\left(1-C_{\omega^M}\right)^{-1}\right|\right|_{L^2 \rightarrow L^2} \leq 
\left(1+||\omega^M||_{L^{\infty}}\right) \frac{1}{1-\left|\left|C_{\omega^M}\right|\right|_{L^2 \rightarrow L^2}^2} \leq c(\varphi), 
\end{equation*}
for $n$ sufficiently large, with $c(\varphi)$ independent of $n$. This proves \eqref{Inversebound}.
The estimate \eqref{keyerrorestimate2} follows similarly from the above estimates. This completes the proof of inequality \eqref{main1}.
Inequalities \eqref{main2} and \eqref{main3} follows directly from \eqref{main1} and the computations 
(with $m = n+1-k$ and $m = n+1-j$)
\begin{equation*}
|||\varphi|||_{0,m} \leq \max_{\gamma \in \{\varphi_+, \varphi_-, \varphi_+^{-1}, \varphi_-^{-1}\}}\, 
\sum_{|l| \geq m}\frac{\nu_l}{e^{|l|\,A(\nu)}}\,|\gamma_l| \leq \frac{1}{e^{m\,A(\nu)}}\,|||\varphi|||_{\nu},
\end{equation*}
and in case $\nu$ increases on $\mathbb{Z}_+$,
\begin{equation*}
|||\varphi|||_{0,m} \leq \max_{\gamma \in \{\varphi_+, \varphi_-, \varphi_+^{-1}, \varphi_-^{-1}\}}\,
\sum_{|l| \geq m}\frac{\nu_l}{\nu_m}\,|\gamma_l| \leq \frac{1}{\nu_m}\,|||\varphi|||_{\nu}.
\end{equation*}
\\ 
We conclude by noting that the above considerations imply the existence of $T_n^{-1}$ for $n$ sufficiently large. Indeed, from the equivalence of 
solvability of RHP's and singular integral equations discussed in Section \ref{IntegrableOp&RHPsection}, it follows from the existence of
$\left(1-C_{\omega^M}\right)^{-1}$ that also $\left(1-C_{\omega}\right)^{-1}$ exists for any factorization 
$v = \left(I - \omega_-\right)^{-1} \left(I + \omega_+\right)$ (with $v$ as in \eqref{vformspec}). So, by the basic relation 
between the integrable operator $K_n$ (as in \eqref{kernelspec}) and the operator $C_{\omega}$ used together with the 
\textit{commutation formula} in \cite{D1} to associate $R_n$ to a RHP, it follows that $\left(1-K_n\right)^{-1}$ exists for $n$ sufficiently
large. Since (according to \eqref{truncatedasintegrable}) the operators $T_n$ and $1-K_n$ agree on $\mathcal{P}_n$, the statement follows. 
Our proof of Theorem \ref{mainresult} is complete.\\
$\Box$

\section{Another look at Baxter's theorem}
\label{Baxtersection}
The following theorem is due to Baxter.
\begin{theorem}
Let $d\mu$ be a non-trivial probability measure on the unit circle and $\nu$ be a strong Beurling weight. Then, 
\begin{equation}
  \label{Baxter}
\sum_{n \in \mathbb{Z}_+} \alpha_n\,z^n \in W_{\nu} \Leftrightarrow d\mu(z) = w(z)\,\frac{|dz|}{2\pi}, w \in W_{\nu}, \min_{z \in \Gamma} w(z) > 0.   
\end{equation}
\end{theorem}

A key element in the proof of inequality \eqref{keyerrorestimate3} lies in the fact that $C_{\omega^M}^2$ (see 
\eqref{Comega}, \eqref{omegaMplusminus}) is a bounded operator in $L^2(\Gamma,|dz|)$ whose norm is small when $n$ is large. 
The same is true for $C_{\omega^M}^2$ as a (bounded) operator in $W_{\nu}$. As we will see, this observation leads to a proof of Theorem
\ref{Baxterextensiontheorem} and thus a new proof of the reverse statement in Baxter's theorem.

\textit{Proof of Theorem \ref{Baxterextensiontheorem}}. Of course, the monic polynomials and hence the Verblunsky coefficients do not change 
if we multiply the weight by a constant: hence we can (and will) assume from the beginning that 
\begin{equation}
  \label{logwassumption}
(\log w)_{0} = 0,
\end{equation} 
without any loss of generality. This will simplify some of the expressions below.

As observed in \cite{D1} the RHP $(\Gamma, v)$, with $v$ as in \eqref{vformspec} (considered in Section \ref{TruncatedToeplitz}), is equivalent
(modulo interchanging $n \leftrightarrow n+1$) to another RHP, namely\\
\quad $\bullet$ $Y_+(z) = Y_-(z)\,\begin{pmatrix} 1 & w/\,z^n\\ 0 & 1\end{pmatrix}$, \quad $z \in \Gamma$,\\
\quad $\bullet$ $Y(z) \begin{pmatrix} z^{-n} & 0\\0 & z^n\end{pmatrix} \rightarrow I$ as $z \rightarrow \infty$.

We shall use the following basic fact: The (1,1)-entry of the (unique) solution of this RHP equals the $n$'th monic OPUC, $Y_{11} = \Phi_n$. 
This RHP, introduced in \cite{BDJ}, is the OPUC analog of the celebrated RHP of Fokas, Its, and Kitaev \cite{FIK} for polynomials orthogonal with 
respect to a weight on the line.

Introduce the successive transformations 
\begin{equation}
Y_1(z) = 
\left\{
\begin{array}{ll}
Y(z)&, |z| < 1,\\
Y(z) \begin{pmatrix} z^{-n} & 0\\0 & z^n\end{pmatrix}&, |z| > 1,   
\end{array}
\right.
\end{equation}
\begin{equation}
Y_2(z) = 
\left\{
\begin{array}{ll}
Y_1(z) \begin{pmatrix} 0 & -1\\1 & 0\end{pmatrix}&, |z| < 1,\\
Y_1(z)&, |z| > 1,   
\end{array}
\right.
\end{equation}
and, with 
\begin{equation}
\delta = \exp \left\{C(\log w)\right\}, \quad \delta_{\pm} = \exp \left\{C_{\pm}(\log w)\right\} \in W_{\nu},
\end{equation}
set
\begin{equation}
Y_3 = Y_2\,\delta^{-\sigma_3}.
\end{equation}
One then easily verifies that (recall \eqref{logwassumption})
\begin{equation}
  \label{PhinY3}
\Phi_n(0) = - \left(Y_3\right)_{12}(0), 
\end{equation}
where $Y_3$ satisfies a normalized RHP $(\Gamma, v_3)$ with jump-matrix
\begin{equation}
v_3 = (I-\omega_-)^{-1}\,(I+\omega_+)
\end{equation} 
and 
\begin{equation}
  \label{omegas}
\omega_- = \begin{pmatrix}0 & 0 \\z^{-n}\,r(z) & 0\end{pmatrix}, \quad
\omega_+ = \begin{pmatrix}0 & -z^n\,r^{-1}(z)\\0 & 0\end{pmatrix}; \quad r = \delta_+^{-1}\,\delta_-^{-1} \in W_{\nu}. 
\end{equation}
By the general theory (recall \eqref{generaltheory}), 
\begin{equation}
  \label{Y3}
Y_3(z) = I + C\left(\mu(\omega_+ +\omega_-)\right)(z), \quad z \in \mathbb{C}\,\backslash\,\Gamma,  
\end{equation}
where
\begin{equation}
  \label{mu}
(1 - C_{\omega})\,\mu = I, \quad \mu \in L^2(\Gamma). 
\end{equation}
It follows from \eqref{PhinY3}, \eqref{omegas} and \eqref{Y3} that 
\begin{equation}
  \label{VerblunskytoRH}
\Phi_n(0)= C\left(\mu_{11}\,z^n\,r^{-1}\right)(0).
\end{equation}
Let us put $\widetilde{\mu}^{(n)} = \mu_{11}$, where we have explicitly indicated the dependence on $n$ in order to avoid confusion in 
the following. It remains to prove that 
\begin{equation}
  \label{remains}
\sum_{n \geq n_0} \nu_n\,\left|C\left(\widetilde{\mu}^{(n)}\,z^n\,r^{-1}\right)(0)\right| < \infty.
\end{equation}
From the first row of \eqref{mu}:
\begin{equation}
\left(\mu_{11}, \mu_{12}\right) = (1,0) + \left(C_+\left(\mu_{12}\,z^{-n}\,r\right), C_-\left(\mu_{11}\,(-z^n\,r^{-1})\right)\right).
\end{equation}
Inserting the equation for $\mu_{12}$ into the equation for $\mu_{11}$ implies the following equation for 
$\widetilde{\mu}^{(n)}$ alone:
\begin{equation}
  \label{mutildeeq}
\widetilde{\mu}^{(n)} = 1 - C_+\left[C_-\left(\widetilde{\mu}^{(n)}\,z^n\,r^{-1}\right)\,z^{-n}\,r\right].
\end{equation}
Clearly,
\begin{equation}
\widetilde{\mu}^{(n)}(z) = \sum_{l \geq 0} \widetilde{\mu}^{(n)}_l\,z^l,
\end{equation}
and we shall write $r(z) = \sum_{k = -\infty}^{\infty} r_k\,z^k$ and $r^{-1}(z) = \sum_{m = -\infty}^{\infty} (r^{-1})_m\,z^m$. 
It follows from \eqref{mutildeeq} that  
\begin{equation}
  \label{mutildecoeffeq}
\widetilde{\mu}^{(n)}_l = \delta_{l,0}\,+ \sum_{p \geq 0, p+m+n<0} (r^{-1})_m\,r_{l-p-m}\,\widetilde{\mu}^{(n)}_p, \quad l \geq 0.
\end{equation}
Let us denote by $W_{\nu}^{\pm}$ the subalgebra of $W_{\nu}$ consisting of functions whose 
negative/non-negative Fourier-coefficients are $0$ and also write $||\cdot||_{\nu^{\pm}} = ||P_{\pm}\cdot||_{\nu}$, where $P_{\pm}$ denotes the 
$L^2$-orthogonal projection onto $H_{\pm}$. Define $\left(A^{(n)}f\right)_l\,$, for $n, l \geq 0$ and $f \in W_{\nu}^{+}$, by 
\begin{equation}
  \label{Andefinition}
\left(A^{(n)}f\right)_l = \sum_{p \geq 0, p+m+n<0} (r^{-1})_m\,r_{l-p-m}\,f_p. 
\end{equation}
With this notation equation \eqref{mutildeeq} takes the form
\begin{equation}
  \label{mutildeintermsofAn}
\widetilde{\mu}^{(n)} = 1 + A^{(n)}\,\widetilde{\mu}^{(n)}.
\end{equation}
Equation \eqref{mutildeintermsofAn} is due essentially to Geronimo and Case (see \cite{Ger-Case}, equations (V.9), (V.10)) and plays an important
role in what follows. The operator $A^{(n)}$ in equation \eqref{mutildeintermsofAn} also appears in \cite{Ger-Case} in a Fredholm determinant 
formula for the Toeplitz determinant  $\det T_n(w)$ (see equation (VII.28)). This formula was rediscovered by Borodin and Okounkov in \cite{B-O} 
and plays an important role in a variety of problems in algebraic combinatorics (see e.g. \cite{BOO}). The operator $A^{(n)}$ is often called the 
Borodin-Okounkov operator.

It is not difficult to establish the following.
\begin{lemma}
  \label{Anlemma}
Let $\nu$ be a Beurling weight and suppose $r \in W_{\nu}$. Then $A^{(n)}$ is a bounded operator on $W_{\nu}^{+}$. Moreover,
$||A^{(n)}||_{W_{\nu}^{+} \rightarrow W_{\nu}^{+}} \rightarrow 0$, as $n \rightarrow \infty$.
\end{lemma}
\textit{Proof.} By submultiplicativity $\nu_l \leq \nu_{l-p-m}\,\nu_p\,\nu_m$, and therefore
\begin{alignat*}{4}
||A^{(n)}f||_{\nu^{+}} &= \sum_{l \geq 0} \nu_l\,\left|\sum_{p \geq 0, p+m+n<0} (r^{-1})_m\,r_{l-p-m}\,f_p\right|\\
&\leq ||r||_{\nu^{+}}\left(\sum_{m<-n} \nu_m\,|(r^{-1})_m|\right)||f||_{\nu^{+}},
\end{alignat*} 
which since $r^{-1} \in W_{\nu}$ proves the claim.\\
$\Box$\\
It follows from \eqref{mutildeintermsofAn} and Lemma \ref{Anlemma} that for $n$ sufficiently large, say $n \geq n_0$,
equation \eqref{mutildeeq} is uniquely solvable and that
\begin{equation}
  \label{nonuniformmuestimate}
||\widetilde{\mu}^{(n)}||_{\nu} \leq c\,||1||_{\nu} \leq c_{\nu}, 
\end{equation}
with a constant $c_{\nu}$ independent of $n$.
We shall need a slightly stronger version of the latter; for $n_0$ sufficiently large
\begin{equation}
  \label{uniformmuestimate}
\sum_{l \geq 0} \nu_l\,\sup_{n \geq n_0}\,|\widetilde{\mu}^{(n)}_l| \leq c_{\nu}.
\end{equation}
To see why this is so, first note from \eqref{mutildecoeffeq} that for $n \geq n_0$
\begin{equation}
  \label{eq1}
|\widetilde{\mu}^{(n)}_l| \leq 
\delta_{l,0}\,+ \sum_{p \geq 0, p+m+n_0<0} |(r^{-1})_m|\,|r_{l-p-m}|\,|\widetilde{\mu}^{(n)}_p|, 
\quad l \geq 0.
\end{equation}
As in the proof of Lemma \ref{Anlemma} we see, that for $n_0$ sufficiently large, the equation
\begin{equation}
  \label{eq2}
s_l = \delta_{l,0}\,+ \sum_{p \geq 0, p+m+n_0<0} |(r^{-1})_m|\,|r_{l-p-m}|\,s_p
\end{equation}
can be (uniquely) solved for $s(z) = \sum_{l \geq 0} s_l\,z^l \in W_{\nu}^{+}$. It suffices to pick $n_0$ so large that the
operator $K : W_{\nu}^{+} \rightarrow W_{\nu}^{+}$ given by 
\begin{equation*}
\left(K f\right)_l = \sum_{p \geq 0, p+m+n_0<0} |(r^{-1})_m|\,|r_{l-p-m}|\,f_p, \quad l \geq 0.
\end{equation*}
has norm less than 1; this is always possible, as in the proof of Lemma \ref{Anlemma}. In the same way that we obtained 
\eqref{nonuniformmuestimate} we see that 
\begin{equation}
  \label{sbound}
||s||_{\nu} \leq c_{\nu}. 
\end{equation}
To prove \eqref{uniformmuestimate} it is therefore enough to show that 
\begin{equation}
  \label{coeffcomparison}
\sup_{n \geq n_0}\,|\widetilde{\mu}^{(n)}_l| \leq s_l, \quad l \geq 0.
\end{equation} 
Denote by $\gamma^{(n)} \in W_{\nu}^{+}$ the element with Fourier coefficients $\gamma^{(n)}_l = |\widetilde{\mu}^{(n)}_l|$, $l \geq 0$. Then we see 
from \eqref{eq1} that 
\begin{equation*}
\gamma^{(n)} + \epsilon^{(n)} = 1 + K\,\gamma^{(n)},
\end{equation*} 
where $\epsilon^{(n)} \in W_{\nu}$ has only \textit{non-negative} Fourier coefficients. 
That is, by \eqref{eq2},
\begin{equation*}
\gamma^{(n)} = (1 - K)^{-1}\left(1 - \epsilon^{(n)}\right) = s - \sum_{j=0}^{\infty} K^{j}\,\epsilon^{(n)}
\end{equation*}
which, since $K$ has non-negative kernel, proves \eqref{coeffcomparison}.\\
Now 
\begin{equation}
  \label{Phin0explicit}
C\left(\widetilde{\mu}^{(n)}\,z^n\,r^{-1}\right)(0) = \int_{\Gamma} \widetilde{\mu}^{(n)}(z)\,z^n\,r^{-1}(z)\frac{dz}{2\pi i z} = 
\sum_{l \geq 0} \widetilde{\mu}^{(n)}_l\,(r^{-1})_{-n-l},
\end{equation}
and we see that to prove \eqref{remains}, it suffices to show that
\begin{equation}
  \label{newremains}
\sum_{n \geq n_0} \sum_{l \geq 0} \nu_n\,|\widetilde{\mu}^{(n)}_l|\,|(r^{-1})_{-n-l}| < \infty.
\end{equation}
But, by \eqref{coeffcomparison}, \eqref{sbound} and the evenness of $\nu$,
\begin{equation*}
\sum_{n \geq n_0} \sum_{l \geq 0} \nu_{n} |\widetilde{\mu}^{(n)}_l|\,|(r^{-1})_{-n-l}| 
\leq \sum_{n \geq n_0} \sum_{l \geq 0} \nu_{n+l} \nu_{-l} s_l |(r^{-1})_{-n-l}| \leq c_{\nu} ||r^{-1}||_{\nu^{-}}.
\end{equation*}
This completes our proof of \eqref{Baxterextension}, and in particular that the RHS of \eqref{Baxter} $\Rightarrow$ LHS of \eqref{Baxter} 
in Baxter's theorem.\\
$\Box$\\
Let us now assume that the Beurling weight $\nu$ is increasing on $\mathbb{Z}_+$. 
Observe first that by \eqref{VerblunskytoRH}, \eqref{mutildecoeffeq} and \eqref{Phin0explicit} we have 
\begin{equation}
  \label{Baxterextensionstart}
\Phi_n(0) = (r^{-1})_{-n} + \sum_{l \geq 0} \left(A^{(n)}\,\widetilde{\mu}^{(n)}\right)_l\,(r^{-1})_{-n-l}.
\end{equation}
By definition \eqref{Andefinition} of $A^{(n)}$ and \eqref{coeffcomparison}, \eqref{sbound}, we have
\begin{alignat*}{4}
&\sum_{n \geq n_0} \nu_n^3\,\left|\sum_{l \geq 0} \left(A^{(n)}\,\widetilde{\mu}^{(n)}\right)_l\,(r^{-1})_{-n-l}\right| =\\
&=\sum_{n \geq n_0} \nu_n^3\,\left|\sum_{p \geq 0, p+m+n<0} (r^{-1})_m\,r_{l-p-m}\,\widetilde{\mu}^{(n)}_p\,(r^{-1})_{-n-l}\right| \\
&\leq \sum_{n \geq n_0} \sum_{l \geq 0} \sum_{p \geq 0} \sum_{m > n+p} 
\nu_{n+l}\,\nu_m\,\nu_{m-p+l}\,|(r^{-1})_{-m}|\,|r_{l-p+m}|\,|(r^{-1})_{-n-l}|\,s_p\\
&\leq c\,||r||_{\nu^{+}}\,||r^{-1}||_{\nu^{-}}^2.
\end{alignat*}
It should be noted that (by first extending the domains of summation) the above sums were carried out by first summing over $n$, then 
over $l$, and finally over $m$ and $p$. This means, by \eqref{Baxterextensionstart}, that
\begin{equation}
  \label{newresult}
\sum_{n \geq n_0} \nu_n^3\,\left|\Phi_n(0) - (r^{-1})_n\right| < \infty.
\end{equation}
It is customary to introduce the \textit{Szeg\"{o} function},
\begin{equation*}
D(z) = \exp\left(\frac{1}{4\pi}\,\int_0^{2\pi} \log w(e^{i\theta})\,\frac{e^{i\theta}+z}{e^{i\theta}-z}\,d\theta\right), \quad z \in 
\mathbb{C}\,\backslash\,\Gamma.
\end{equation*} 
Note that $D(z)$ and $\delta(z)$ are in general proportional, and that in case $(\log w)_0 = 0$ (see above) they are equal.
Following Simon we also introduce the function
\begin{equation*}
S(z) = -\sum_{n=1}^{\infty} \alpha_{n-1}\,z^n,
\end{equation*}
where $\alpha_{n-1} \equiv - \overline{\Phi_n(0)}$ for $n \geq n_0$ and $\alpha_{n-1} \equiv 0$ for $n < n_0$.
We shall use the notation $D_i$ resp. $D_e$ for the restriction of $D$ to the interior resp. exterior of the unit circle, as well as for the 
analytical continuations of these functions across the unit circle, should they exist.
Now, $(r^{-1})_{-n} = \overline{\left((\overline{r}^{-1})_n\right)}$ and $\overline{r}^{-1} = \overline{D}_i\,\overline{D}_e$. 
Also, if $w$ is positive, then $D_e(z) = 1/\overline{D_i(1/\overline{z})}$, $|z| > 1$.
Equation \eqref{newresult} therefore implies the following result.
\begin{theorem}
  \label{Simon-improvement}
Let $\nu$ be a Beurling weight which increases on $\mathbb{Z}_+$ and $d\mu(z) = w(z)\,|dz|$ be a measure 
on the unit circle. Suppose that $w \in W_{\nu}$, $w \not= 0$ on $\mathcal{R}_{\nu}$, $wind(w, 0) = 0$. Then,
\begin{equation}
\overline{D}_i\,\overline{D}_e - S \in W_{\nu^3}^{+}.
\end{equation}
In particular, for $w$ positive, we obtain 
\begin{equation}
  \label{Bornapprox}
\frac{\overline{D}_i}{D_i} - S \in W_{\nu^3}^{+}.
\end{equation}
\end{theorem}
This theorem should be viewed as a refinement of the reverse implication in Baxter's theorem: not only is $S \in W_{\nu}$, but 
$S = \overline{D}_i\,\overline{D}_e$ up to three orders of smoothness. Alternatively, from a physical point of view we can regard 
$\overline{D}_i\,\overline{D}_e$ as the principal object of study: indeed for real weights, $r = \overline{r}^{-1} = \frac{\overline{D}_i}{D_i}$
is the reflection coefficient for the system at hand and $S$ is the leading Born approximation (see \cite{Sim2}, \cite{Sim3}). Thus, 
\eqref{Bornapprox} is an estimate of how the Born approximation deviates from $r$.

It is a well-known theorem of Nevai and Totik (\cite{NT}) that for real $d\mu$, $\limsup_{n \rightarrow \infty} |\alpha_n|^{1/n} = R^{-1} < 1$ 
if and only if $d\mu$ obeys the Szeg\"{o} condition, $d\mu_s = 0$ and $D_i^{-1}$ has an analytic extension to 
$\left\{z \in \mathbb{C} : |z| < R\right\}$. Theorem \ref{Simon-improvement} therefore has the following corollary. 
\begin{corollary}
  \label{Analyticity-corollary}
Let $d\mu$ be a positive measure on $\Gamma$. Suppose that 
\begin{equation}
\limsup_{n \rightarrow \infty} |\alpha_n|^{1/n} = R^{-1} < 1,
\end{equation}
so that $D_i^{-1}$ and $S$ are analytic in $\left\{z \in \mathbb{C} : |z| < R\right\}$. Then, for some $\delta > 0$, the function
$\overline{D_i(\frac{1}{\overline{z}})}/D_i(z)-S(z)$ is analytic in $\left\{z \in \mathbb{C} : 1-\delta < |z| < R^3\right\}$.
\end{corollary} 
\textit{Proof.} It follows from the result of Nevai and Totik that $\frac{1}{w} = \frac{D_e}{D_i}$ is analytic, and in particular that 
$w$ cannot vanish, in the set $\left\{z \in \mathbb{C} : 1/R < |z| < R \right\}$. In addition, as $w > 0$ on $\Gamma$, $wind(w,0) = 0$. 
We may then, for any $\epsilon > 0$, apply Theorem \ref{Simon-improvement} to the Beurling weight defined by 
$v_{n} = \left(R\,(1-\epsilon)\right)^{|n|}$ for $n \in \mathbb{Z}$. This proves analyticity in $\{1 < |z| < R^3\}$. The analyticity in
$\{1-\delta < |z| < R^3\}$ follows from the fact that $D_i$ is meromorphic in $|z| < R$, but has no poles on $\Gamma$.\\
$\Box$\\   
In \cite{Sim1} Simon proved Corollary \ref{Analyticity-corollary} with $R^3$ replaced by $R^2$, see Theorem 7.2.1. Motivated by Corollary
\ref{Analyticity-corollary} above, Simon \cite{Sim3} has now given an independent proof of the result.
\section{Some examples}
\label{Examplesection}
We thank Barry Simon for drawing our attention to the following examples from \cite{Sim1}, which illustrate the sharpness of Corollary
\ref{Analyticity-corollary} (see also \cite{Sim3}).
\\
\textit{Example 1 (Single nontrivial moment).} Consider the weight $w(e^{i\theta}) = 1 - a \cos{\theta}$, $0 < a < 1$, 
having a single nontrivial moment. Introduce the auxiliary parameters
\begin{equation}
\mu_{\pm} = a^{-1} \pm \sqrt{a^{-2}-1}. 
\end{equation}
Note that $\mu_+\,\mu_- = 1, 0 < \mu_- < 1$. By computation one finds that
\begin{equation}
D_i(z) = \sqrt{\frac{a}{2\,\mu_{-}}}\,\left(1-\frac{z}{\mu_{+}}\right),
\end{equation}
and so $D_i^{-1}$ has a simple pole at $z=\mu_+$. Also, 
\begin{alignat}{4}
\alpha_{n} &= -\frac{\mu_+-\mu_-}{\mu_+^{n+2}-\mu_-^{n+2}} = -(\mu_+-\mu_-)\,\mu_+^{-n-2}\,\left(1-\mu_+^{-(2n+4)}\right)^{-1}\\
\notag &= -(\mu_+-\mu_-)\,\sum_{j=1}^{\infty} (\mu_+^{-n-2})^{2j-1}, 
\end{alignat}
so that $S$ has simple poles at $z_{j} = \mu_+^{2j-1}$, $j \in \mathbb{N}$. The statement in Corollary \ref{Analyticity-corollary} is easily
verified by noting that $Res(\overline{D_i(1/\overline{z})}/D_i(z), z = \mu_+) = Res(S, z = \mu_+) = -(\mu_+-\mu_-)$.\\
\textit{Example 2 (Rogers-Szeg\"{o} polynomials).} Let $0 < q < 1$ and consider the weight with Verblunsky coefficients 
\begin{equation}
\alpha_n = (-1)^n\,q^{(n+1)/2}, \quad n \geq 0.
\end{equation}
Then, 
\begin{equation}
D_i(z) = \Pi_{j=0}^{\infty}\,(1-q^{j+1})^{1/2}\,(1+q^{j+1/2}\,z)
\end{equation}
so that $D_i^{-1}$ has simple poles at $z_j = -q^{-j-1/2}$, $j \geq 0$. On the other hand,
\begin{equation}
S(z) = -\sum_{n=1}^{\infty}\,(-1)^n\,q^{n/2}\,z^n = -\frac{q^{1/2}\,z}{1+q^{1/2}\,z}
\end{equation}
has a simple pole at $z = -q^{-1/2}$. The statement in Corollary \ref{Analyticity-corollary} follows from
$Res(\overline{D_i(1/\overline{z})}/D_i(z), z = -q^{-1/2}) = Res(S, z = -q^{-1/2}) = q^{-1/2}$.\\
\section{The inverse statement in a theorem of Golinskii-Ibragimov}
\label{Ibragimovsection}
Let us denote by $H^{1/2}$ the Sobolev space of functions $f \in L^{2}(\Gamma)$ with $\sum_{l \in \mathbb{Z}} |l|\,|f_l|^2 < \infty$, equipped
with the norm $||f||_{1/2} = (\sum_{l \in \mathbb{Z}} (1+|l|)\,|f_l|^2)^{1/2}$. Let $H^{1/2}_{\mathbb{R}}$ denote the class of real-valued 
functions in $H^{1/2}$. The following theorem is implied by the Ibragimov/Golinskii-Ibragimov version of the Strong Szeg\"{o} Limit Theorem \cite{Sim1}.
\begin{theorem}
  \label{Ibragimovtheorem}
Let $d\mu$ be a non-trivial probability on the unit circle. Then,
\begin{equation}
  \label{Ibragimoveq}
\sum_{n \in \mathbb{Z}_+} n\,|\alpha_n|^2 < \infty \Leftrightarrow d\mu = w\,\frac{|dz|}{2\pi} \mbox{ and } \log w \in H^{1/2}_{\mathbb{R}}.
\end{equation}  
\end{theorem} 
Just as Riemann-Hilbert techniques provide a direct proof of the $\mathcal{I}$-part of Baxter's theorem, they can also be used to proof
that the RHS of \eqref{Ibragimoveq} $\Rightarrow$ LHS of \eqref{Ibragimoveq}. This is the goal of this section.\\
We will need the following proposition (see \cite{Sim1}, Prop. 6.2.6).
\begin{proposition}
  \label{Deift-Killip}
For $f \in H^{1/2}_{\mathbb{R}}$, let 
\begin{equation}
I(f) = -\sum_{k > 0} f_k\,z^k + \sum_{k < 0} f_k\,z^k
\end{equation}
and 
\begin{equation}
B(f) = \exp(I(f)).
\end{equation}
Then $B$ maps $H^{1/2}_{\mathbb{R}}$ continuously into $H^{1/2}$.
\end{proposition}
It follows immediately from the above that if $\log w \in H^{1/2}_{\mathbb{R}}$, then $r = \frac{\overline{D}_i}{D_i} = B(\log w) \in H^{1/2}$.
Next observe that for real measures $d\mu$, $\overline{r}_m = r_{-m}$, and hence \eqref{Andefinition} takes the form
\begin{equation}
  \label{Borodin-Okounkovoperator}
(A^{(n)}f)_l = \sum_{p \geq 0} \left(\sum_{m > n}\,r_{l+m}\,\overline{r_{m+p}}\right) f_p, \quad l \geq 0.
\end{equation}
Previously we regarded $A^{(n)}$ as an operator in $W_{\nu}$. However, $A^{(n)}$ can also be regarded as a trace class 
(and in particular bounded), positive, self-adjoint operator on $l^2_+ \equiv l^2(\mathbb{Z}_+) \cong H_+$. Indeed, $A^{(n)}$ has the 
form $R\,\chi_n\,R^{*}$ where $R$ is the Hilbert-Schmidt operator on $l^2_+$ with kernel $R_{i,j} = r_{i+j}$, $i,j \geq 0$,
\begin{equation}
  \label{RHilbSchmidt}
||R||_{\mathcal{I}_2(l^2_+)}^2 = \sum_{i,j \geq 0} |r_{i+j}|^2 = \sum_{i \geq 0} (1+i)\,|r_i|^2 \leq ||r||^2_{1/2}
\end{equation}
and $\chi_n$ denotes multiplication by the characteristic function of the set $\{m > n\}$. It follows that $A^{(n)}$ is trace class in $l^2_+$ with
\begin{equation}
  \label{Annorm}
||A^{(n)}||_{l^2_+ \rightarrow l^2_+} \leq ||A^{(n)}||_{\mathcal{I}_1(l^2_+)} = \sum_{l \geq 0} \sum_{m > n} |r_{l+m}|^2 \leq 
\sum_{m > n} (1+m)\,|r_m|^2.
\end{equation} 
From \eqref{VerblunskytoRH}, \eqref{Phin0explicit} 
\begin{equation}
  \label{Almostscalarprod}
\alpha_{n-1} = - \overline{\Phi_n(0)} = -\sum_{l \geq 0} \overline{\widetilde{\mu}^{(n)}_l}\,r_{n+l}. 
\end{equation}
Here $\widetilde{\mu}^{(n)} = (\widetilde{\mu}^{(n)}_l)_{l \geq 0}$ solves the equation \eqref{mutildeintermsofAn} in $W_{\nu}$. However,
by \eqref{Annorm} equation \eqref{mutildeintermsofAn} is also uniquely solvable in $l^2_+$ for $n$ sufficiently large. 
As $W_{\nu} \hookrightarrow l^2_+$, it follows that we may regard $\widetilde{\mu}^{(n)}$ as the (unique) solution of 
\eqref{mutildeintermsofAn} in $l^2_+$. 
But $r^{(n)} = (r_{n+l})_{l \geq 0}$ is also in $W_{\nu} \hookrightarrow l^2_+$ and hence we may write \eqref{Almostscalarprod} in the form
\begin{equation}
  \label{Verblunskyformula}
\alpha_{n-1} = - \left(r^{(n)},\frac{1}{1-A^{(n)}}\,e_0\right)_{l^2_+},
\end{equation}
where $e_0 = (1, 0, 0,...)^{T}$ and the inverse of $1 - A^{(n)}$ is taken in $l^2_+$. 

Equation \eqref{Verblunskyformula} is derived in the case $w \in W_{\nu}$, but as we now show, it remains true for $w$ with 
$\log w \in H^{1/2}_{\mathbb{R}}$. Note first that for $f \equiv \log w \in H^{1/2}_{\mathbb{R}}$, $w \in L^p(\Gamma)$ for all $1 \leq p < \infty$ by 
(the proof of) Lemma 6.1.4 in \cite{Sim1}. Set $f^{(N)} = \sum_{-N}^{N} f_j\,z^j$ and $w^{(N)} = e^{f^{(N)}} \in W_{\nu}$ for any Beurling 
weight $\nu$. If $D_i^{(N)}$ denotes the Szeg\"{o} function for $w^{(N)}$, then $r^{(N)} = \frac{\overline{D_i^{(N)}}}{D_i^{(N)}} = B(f^{(N)})$.
Let $A^{(N,n)}$ denote the Borodin-Okounkov operator \eqref{Borodin-Okounkovoperator} with $r$ replaced by $r^{(N)}$. By \eqref{Annorm}
\begin{alignat*}{4}
||A^{(N,n)}||_{l^2_+ \rightarrow l^2_+}^{1/2} &\leq \bigg(\sum_{m > n} (1+m)|r^{(N)}_m|^2\bigg)^{1/2} \\
\notag
&\leq \bigg(\sum_{m \geq 0} (1+m)|r^{(N)}_m - r_m|^2\bigg)^{1/2} + \bigg(\sum_{m > n} (1+m)|r_m|^2\bigg)^{1/2}. 
\end{alignat*}
The first term on the right converges to zero as $N \rightarrow \infty$ by Proposition \ref{Deift-Killip}, and the second term can be made
small uniformly for $n$ large. Thus, for any fixed $\rho_0 < 1$, there exists $N_0, n_0$ such that 
\begin{equation}
  \label{AnNuniformbound}
||A^{(N,n)}||_{l^2_+ \rightarrow l^2_+} < \rho_0^2
\end{equation}
if $N \geq N_0$ and $n \geq n_0$. Hence for all $N \geq N_0$ and $n \geq n_0$ we have by \eqref{Verblunskyformula}
\begin{equation}
  \label{Verblunskyformula-approx}
\alpha_{n-1}^{(N)} = - \left(r^{(N,n)},\frac{1}{1-A^{(N,n)}}\,e_0\right)_{l^2_+},
\end{equation}
where $r^{(N,n)} = (r^{(N)}_{n+l})_{l \geq 0}$ and $\alpha_{n-1}^{(N)}$ is the $(n-1)^{st}$ Verblunsky coefficient for $w^{(N)}$. But for fixed $n$, 
a simple computation shows that as $N \rightarrow \infty$, $r^{(N,n)} \rightarrow r^{(n)}$ in $H^{1/2} \hookrightarrow l^2_+$, 
and in addition, by \eqref{RHilbSchmidt}, $A^{(N,n)} \rightarrow A^{(n)}$ in $\mathcal{I}_1(l^2_+) \subset \mathcal{L}(l^2_+)$, the bounded 
operators on $l^2_+$. Finally, using \eqref{AnNuniformbound}, we see that for all $n \geq n_0$ the RHS of \eqref{Verblunskyformula-approx} 
converges to the RHS of \eqref{Verblunskyformula}. But as $N \rightarrow \infty$ the LHS of \eqref{Verblunskyformula-approx} converges to the LHS of 
\eqref{Verblunskyformula} by Lemma 6.1.4 (b) in \cite{Sim1}. This establishes \eqref{Verblunskyformula} for $w$ with $\log w \in H^{1/2}_{\mathbb{R}}$
and $n \geq n_0$.

\textit{Remark.} The reader may ask why we do not prove \eqref{Verblunskyformula} directly from the RHP in Section \ref{Baxtersection} with 
weight $w$, $\log w \in H^{1/2}_{\mathbb{R}}$, rather than proceeding by approximation as above. However, we only know that 
$w \in L^p(\Gamma)$ for $1 \leq p < \infty$, not in $L^{\infty}(\Gamma)$. Thus the RHP is non-standard and requires special (BMO) considerations, which
we can, and do, avoid.\\

We will now show that $\sum_{n \in \mathbb{Z}_+} n\,|\alpha_n|^2 < \infty$. Note first from \eqref{Annorm}, that for $n \geq n_0$ 
\begin{equation}
  \label{Anbound}
||A^{(n)}||_{l^2_+ \rightarrow l^2_+} \leq \rho^2 
\end{equation}
where
\begin{equation}
\rho \equiv \bigg(\sum_{n \geq n_0} (n+1)\,|r_n|^2\bigg)^{1/2} < 1.
\end{equation}
Secondly, using formula \eqref{Verblunskyformula} we obtain
\begin{equation*}
\alpha_{n-1} = -r_n - \big(r^{(n)}, (1-A^{(n)})^{-1}\,A^{(n)}\,e_0\big)_{l^2_+},
\end{equation*}
and therefore
\begin{alignat}{4}
  \label{Ibragimov1}
\bigg(\sum_{n \geq n_0} &n\,|\alpha_{n-1}|^2\bigg)^{1/2} \leq 
\rho + \bigg(\sum_{n \geq n_0} n\,\big|\big(r^{(n)}, (1-A^{(n)})^{-1}\,A^{(n)}\,e_0\big)_{l^2_+}\big|^2\bigg)^{1/2} \\
\notag
&\leq \rho + \bigg(\sum_{n \geq n_0} n\,||r^{(n)}||_{l^2_+}^2\,||(1-A^{(n)})^{-1}\,A^{(n)}\,e_0||_{l^2_+}^2\bigg)^{1/2} \\
\notag
&\leq \rho + \bigg(\sup_{n \geq n_0} n\,||(1-A^{(n)})^{-1}\,A^{(n)}\,e_0||_{l^2_+}^2\bigg)^{1/2} \cdot 
\bigg(\sum_{n \geq n_0} ||r^{(n)}||_{l^2_+}^2\bigg)^{1/2}.
\end{alignat}
Obviously,
\begin{equation}
  \label{Ibragimov2}
\sum_{n \geq n_0} \big|\big|r^{(n)}\big|\big|_{l^2_+}^2 = \sum_{n \geq n_0} \sum_{j \geq 0} |r_{n+j}|^2 \leq 
\sum_{j \geq n_0} (j+1)\,|r_j|^2 = \rho^2.
\end{equation}
Furthermore, by \eqref{Anbound}, for any $n \geq n_0$
\begin{alignat}{4}
\notag
||(1-A^{(n)})^{-1}\,A^{(n)}\,e_0||_{l^2_+}^2 &\leq ||(1-A^{(n)})^{-1}||_{l^2_+ \rightarrow l^2_+}^2\,||A^{(n)}\,e_0||_{l^2_+}^2 \\
\label{Ibragimov3}
&\leq \frac{1}{(1-\rho^2)^2}\,||A^{(n)}\,e_0||_{l^2_+}^2
\end{alignat}
and also
\begin{alignat}{4}
\notag
n\,\big|\big|A^{(n)}\,e_0\big|\big|_{l^2_+}^2 &= n\,\sum_{l \geq 0} \bigg|\sum_{m > n} r_{l+m}\,\overline{r}_m\bigg|^2 \leq 
n\,\sum_{l \geq 0} \bigg(\sum_{m > n} |r_{l+m}|^2\bigg) \bigg(\sum_{m > n} |r_m|^2\bigg) \\
\label{Ibragimov4}
&\leq n\,\rho^2\,\sum_{m > n} |r_m|^2 \leq \rho^2\,\sum_{m > n} m\,|r_m|^2 \leq \rho^4.
\end{alignat}
It follows from \eqref{Ibragimov3} and \eqref{Ibragimov4}, that
\begin{equation}
  \label{Ibragimov5}
\sup_{n \geq n_0} n\,||(1-A^{(n)})^{-1}\,A^{(n)}\,e_0||_{l^2_+}^2 \leq \frac{1}{(1-\rho^2)^2}\,\rho^4.
\end{equation}
Combining \eqref{Ibragimov1}, \eqref{Ibragimov2} and \eqref{Ibragimov5}, it follows that
\begin{equation}
\bigg(\sum_{n \geq n_0} n\,|\alpha_{n-1}|^2\bigg)^{1/2} \leq \frac{\rho}{1-\rho^2}.
\end{equation}
This completes the proof that the RHS of \eqref{Ibragimoveq} $\Rightarrow$ LHS of \eqref{Ibragimoveq}.\\
$\Box$

\textit{Acknowledgements.} The work of the first author was supported in part by the NSF Grant DMS-0296084.
The second author would like to express his gratitude to the Wenner-Gren Foundations for  
their financial support. He is also grateful for the hospitality and stimulating environment provided by the Courant Institute.
Both authors would like to thank Barry Simon for fruitful discussions. 


\begin{thebibliography}{99}
\bibitem{BDJ} J. Baik, P. Deift, K. Johansson, \textit{On the distribution of the length of the longest increasing subsequence of random permutations},
J. Amer. Math. Soc. \textbf{12} (1999), no 4, 1119-1178.
\bibitem{B-O} A.M. Borodin, A. Okounkov, \textit{A Fredholm determinant formula for Toeplitz determinants}, Integral Equations and Operator Theory
\textbf{37} (2000), 386-396.
\bibitem{BOO} A. Borodin, A. Okounkov, G. Olshanski, \textit{Asymptotics of Plancherel measures for symmetric groups}, J. Amer. Math. Soc.
\textbf{13} (2000), 491-515.
\bibitem{BottSilb} A. B\"{o}ttcher, B. Silbermann, \textit{Introduction to large truncated Toeplitz matrices}, Springer-Verlag (1999).
\bibitem{CG} K. Clancey, I. Gohberg, \textit{Factorization of matrix functions and singular integral operators}, Operator Theory, vol. 3, 
Birkh\"{a}user, Basel, 1981. 
\bibitem{CMV} M.J. Cantero, L. Moral, L. Vel\'{a}zquez, \textit{Five-diagonal matrices and zeros of orthogonal polynomials on the unit circle},
Linear Algebra  Appl. \textbf{362} (2003), 29-56. 
\bibitem{D1} P. A. Deift, \textit{Integrable Operators}, M. Sh. Birman's 70th anniversary collection
(V. Buslaev, M. Solomjak, D. Yafaev, eds.),  Amer. Math. Soc. Transl. ser. 2, \textbf{159}, Amer. Math. Soc., Providence 1999. 
\bibitem{D2} P. A. Deift, \textit{Orthogonal Polynomials and Random Matrices: A Riemann-Hilbert Approach}, Courant Lecture Notes \textbf{3} (1999). 
\bibitem{Doug} R. G. Douglas, \textit{Banach Algebra Techniques in Operator Theory}, Academic Press, New York 1972. 
\bibitem{DZ1} P. A. Deift, X. Zhou, \textit{A steepest descent method for oscillatory Riemann-Hilbert problems. Asymptotics for the MKdV
equations}, Ann. of Math. \textbf{137} (1993), 295-368.
\bibitem{FIK} A. S. Fokas, A. R. Its, A. V. Kitaev, \textit{Discrete Painleve' equations and their appearance in quantum gravity}, Comm. Math. Phys.
\textbf{142} (1991), 313-344. 
\bibitem{Ger-Case} K. M. Case, J.S. Geronimo, \textit{Scattering theory and polynomials orthogonal on the unit circle}, J. Math. Phys. 
\textbf{20}(2) (1979)
\bibitem{IIKS}  A. R. Its, A. G. Izergin, V. E. Korepin and N. A. Slavnov, \textit{Differential equations for quantum correlation functions}, Int. J.
Mod. Phys. \textbf{B4} (1990), 1003-1037; \textit{The quantum correlation function as the $\tau$ function of classical differential equations}, 
Important developments in soliton theory, (A. S. Fokas and V. E. Zakharov, eds.), Springer-Verlag, Berlin (1993), 407-417.  
\bibitem{Joh} K. Johansson, \textit{Discrete Polynuclear Growth and Determinantal Processes}, Comm. Math. Phys. \textit{242} (2003), 277-329.
\bibitem{NT} P. Nevai, V. Totik, \textit{Orthogonal polynomials and their zeros}, Acta Sci. Math. \textbf{53} (1989), 99-104.
\bibitem{Sak} L.A. Sakhnovich, \textit{Operators similar to unitary operators}, Functional Anal. and Appl. \textbf{2} No. 1 (1968), 48-60.
\bibitem{Sim1} B. Simon, \textit{Orthogonal Polynomials on the Unit Circle, Part 1: Classical Theory}, AMS Colloquium Series, Vol. \textbf{54}, 
Amer. Math. Soc., Providence, RI, 2005. 
\bibitem{Sim2} B. Simon, \textit{Orthogonal Polynomials on the Unit Circle, Part 2: Spectral Theory}, AMS Colloquium Series, Vol. \textbf{54},
Amer. Math. Soc., Providence, RI, 2005. 
\bibitem{Sim3} B. Simon, \textit{Meromorphic Szeg\"{o} functions and asymptotic series for Verblunsky coefficients}, in preparation.
\bibitem{Sz} G. Szeg\"{o}, \textit{Orthogonal polynomials}, AMS Colloquium Series, Vol. \textbf{23}, Amer. Math. Soc., 
Providence, RI, 1939; 3rd edition 1967.
\end{thebibliography}
\end{document}